\documentclass[11pt,a4paper]{article}
\usepackage[T1]{fontenc}
\usepackage{lmodern,microtype}
\usepackage[margin=28mm]{geometry}
\usepackage{amsmath,amssymb,amsthm,mathtools,bm}
\usepackage{mathrsfs}
\usepackage{booktabs,array}
\usepackage{etoolbox}
\usepackage{xcolor}
\usepackage{tikz}
\usepackage{pgfplots}
\pgfplotsset{compat=1.18}
\usetikzlibrary{arrows.meta,calc}
\usepackage[font=small,labelfont=bf]{caption}
\definecolor{metalblue}{RGB}{34,89,123}
\definecolor{metalwash}{RGB}{226,237,245}
\definecolor{screenorange}{RGB}{190,89,29}
\definecolor{screentint}{RGB}{248,229,208}
\definecolor{fieldteal}{RGB}{0,112,111}
\tikzset{figarrow/.style={-{Latex[length=2mm]},semithick},
         figdim/.style={{Latex[length=1.5mm]}-{Latex[length=1.5mm]},thin},
         figboundary/.style={draw=metalblue,thick}}
\usepackage[colorlinks=true,linkcolor=blue!45!black,citecolor=blue!45!black,urlcolor=blue!45!black]{hyperref}
\hypersetup{pdftitle={Quantitative surface approximations for Maxwell scattering in nonlocal plasmonics},
pdfauthor={Habib Ammari, Hyundae Lee, Matias Ruiz, Sanghyeon Yu},
pdfkeywords={nonlocal plasmonics, Maxwell scattering, hydrodynamic Drude model, GNOR, effective transmission conditions}}
\numberwithin{equation}{section}
\newtheorem{theorem}{Theorem}[section]
\newtheorem{proposition}[theorem]{Proposition}
\newtheorem{lemma}[theorem]{Lemma}
\newtheorem{corollary}[theorem]{Corollary}
\theoremstyle{definition}

\theoremstyle{remark}
\newtheorem{remark}[theorem]{Remark}
\newcommand{\R}{\mathbb R}
\newcommand{\C}{\mathbb C}
\newcommand{\Sph}{\mathbb S^2}
\newcommand{\ii}{\mathrm i}
\newcommand{\dd}{\,\mathrm d}
\newcommand{\Id}{\mathrm{Id}}
\newcommand{\eps}{\varepsilon}
\newcommand{\jump}[1]{\mathopen{[\![}#1\mathclose{]\!]}}
\newcommand{\norm}[1]{\lVert #1\rVert}

\newcommand{\E}{\mathbf E}
\newcommand{\HH}{\mathbf H}
\newcommand{\J}{\mathbf j}
\newcommand{\ET}{\mathbf E^{\mathrm{T}}}
\newcommand{\nn}{\mathbf n}
\newcommand{\GG}{\mathcal G_\kappa}
\newcommand{\NN}{\mathcal N_\eta}
\newcommand{\NNtwo}{\widehat{\mathcal N}^{(2)}_\eta}
\newcommand{\TT}{\mathcal T_z}
\newcommand{\LL}{\mathcal L_z}
\newcommand{\FF}{\mathbf F}
\newcommand{\inc}{\mathrm{inc}}
\newcommand{\scat}{\mathrm{sc}}
\newcommand{\eff}{\mathrm{eff}}
\newcommand{\gnor}{\mathrm{gnor}}
\newcommand{\diff}{\zeta}
\newcommand{\sg}{s_{\diff}}
\newcommand{\Qg}{Q_{\diff}}
\newcommand{\varg}{\vartheta_{\diff}}
\DeclareMathOperator{\diver}{div}
\DeclareMathOperator{\curl}{curl}

\DeclareMathOperator{\diag}{diag}
\DeclareMathOperator{\ran}{ran}
\allowdisplaybreaks[2]
\title{Quantitative surface approximations\\
for Maxwell scattering in nonlocal plasmonics}
\author{Habib Ammari\textsuperscript{1,2}\quad Hyundae Lee\textsuperscript{3}\quad
Matias Ruiz\textsuperscript{4}\quad Sanghyeon Yu\textsuperscript{5}}
\date{}
\makeatletter
\renewcommand{\maketitle}{\begin{center}
  {\LARGE\@title\par}\vspace{0.8em}
  {\normalsize\@author\par}\vspace{0.8em}
  {\footnotesize
  \textsuperscript{1}Department of Mathematics, ETH Z\"urich,
  R\"amistrasse 101, 8092 Z\"urich, Switzerland\\
  \textsuperscript{2}Hong Kong Institute for Advanced Study,
  City University of Hong Kong, Kowloon Tong, Hong Kong\\
  \href{mailto:habib.ammari@math.ethz.ch}{\texttt{habib.ammari@math.ethz.ch}}\\[3pt]
  \textsuperscript{3}Department of Mathematics, Inha University, Incheon, Korea\\
  \href{mailto:hdlee@inha.ac.kr}{\texttt{hdlee@inha.ac.kr}}\\[3pt]
  \textsuperscript{4}School of Computing and Mathematical Sciences,
  University of Leicester, Leicester, UK\\
  \href{mailto:mr447@leicester.ac.uk}{\texttt{mr447@leicester.ac.uk}}\\[3pt]
  \textsuperscript{5}Department of Mathematics, Korea University, Seoul, Korea\\
  \href{mailto:sanghyeon_yu@korea.ac.kr}{\texttt{sanghyeon\_yu@korea.ac.kr}}\par}
  {\small\@date\par}
\end{center}\vspace{0.3em}}
\makeatother

\begin{document}
\maketitle
\begin{abstract}
We establish quantitative surface approximations for Maxwell scattering
by a fixed collection of smooth, separated metallic particles described
by the hard-wall generalized nonlocal optical response (GNOR) model.
As the ratio $\eta$ of the reference screening length to the particle
size tends to zero, regularized first-order and second-order
transmission conditions approximate the electric far-field operator
with errors $O(\eta^2)$ and $O(\eta^3)$, respectively. The estimates
retain retardation and are uniform on compact subplasma frequency
intervals with fixed positive damping and for a bounded dimensionless
diffusion-to-pressure ratio. The proof combines exact elimination of the
longitudinal field, a uniform weak bound for its boundary density, and
smoothing by exterior observations. The same error orders hold for
fields observed at a fixed positive distance from the particles.
We also justify thin-coating representations on bounded tangential
wave-number ranges for a planar interface and on fixed finite spherical
angular subspaces. For each fixed quasistatic spherical multipole with a
simple local pole, the first- and second-order models approximate the pole
and its residue with errors $O(\eta^2)$ and $O(\eta^3)$, respectively.
\end{abstract}




\section{Introduction}\label{sec:introduction}

Plasmonic resonances arise from collective oscillations of conduction
electrons in metallic particles and produce strong electromagnetic field
enhancement. For particles described by local permittivity, their
mathematical analysis is well developed. Layer potential techniques
relate the resonances to the spectrum of the Neumann--Poincar\'e
operator and therefore to the geometry of the particles
\cite{AmmariDengMillien,AmmariScalar,AmmariMaxwell,AmmariBook}.
Many of the most striking effects are multiscale. In nearly touching
particles, for instance, the field concentrates in the gap, on a length
much smaller than the particle size
\cite{YuAmmariNanospheres,YuAmmariHybridization}. When the scales
are separated, asymptotic reductions and effective models become central to
the analysis of resonant structures \cite{AmmariCBMS}.

Electronic screening introduces a further small scale, internal to each
particle. Nonlocal effects appear when the electric field varies on
this scale, as it does near small particles, sharp curvature, or narrow
gaps. A local permittivity is then no longer sufficient and the
constitutive law must account for spatial dispersion. A convenient
semiclassical description is the generalized nonlocal optical response
(GNOR) model \cite{GNOR}. It treats the conduction electrons as
charged fluids whose current responds to the electric field, a
pressure gradient, and the diffusion of an induced charge. The pressure
confines the induced charge to a screening layer at the metal surface.
Diffusion damps this charge and accounts for size-dependent broadening
\cite{StamatopoulouTserkezis,SvendsenGNOR}. A dimensionless parameter
$\diff\geq0$ measures the diffusion relative to pressure, and $\diff=0$
gives the widely used hydrodynamic Drude model
\cite{RazaReview,MortensenReview}. The hard-wall form of the model
forbids electrons to leave the metal. Self-consistent and
quantum-corrected models, which describe inhomogeneous equilibrium
densities, spill-out, and tunnelling \cite{Toscano,Esteban}, lie beyond
this description.

The aim of this paper is to justify a surface approximation of Maxwell
scattering by metallic particles whose electronic response is described
by the hard-wall GNOR model. The approximation replaces the screening
layer by a transmission condition on the particle surface. Let $a$ be a
reference particle size, $\ell$ the reference screening length, and
$\eta=\ell/a$. We keep the smooth particles and their positive separations fixed as
$\eta\to0$. Closing gaps introduce a second small length
\cite{SchnitzerPRB}. We fix the positive damping coefficient, take
$0\leq\diff\leq\diff_*$ with $\diff_*<\infty$, and consider frequencies
$z=\omega/\omega_p$ in a compact interval $I\Subset(0,1)$, where
$\omega_p$ is the plasma frequency. Let $\mathbb S^2$ be the unit
sphere and let $L_t^2(\mathbb S^2)$ be the space of square-integrable
tangential vector fields on it. For an incident electric Herglotz wave
with density $\mathbf a\in L_t^2(\mathbb S^2)$, denote by
$\FF_\eta^{\gnor}$ and $\FF_\eta^{\eff}$ the electric far fields of the
GNOR and the effective problems, respectively.
Theorem~\ref{thm:farfield} proves
\[
\sup_{z\in I,\ 0\leq\diff\leq\diff_*}
\norm{\FF_\eta^{\gnor}-\FF_\eta^{\eff}}_{L^2(\mathbb S^2)}
\leq C\eta^2\norm{\mathbf a}_{L^2(\mathbb S^2)}.
\]
The effective problem consists of local Maxwell equations coupled by a
regularized first-order transmission condition. The estimate retains
retardation and imposes no finite-mode restriction on the incident field.
Theorem~\ref{prop:secondordersurface} constructs a regularized
second-order surface model whose far-field error is $O(\eta^3)$ under
the same assumptions. Diffusion adds a nonnegative term to the
current dissipation and modifies the longitudinal screening parameter,
which remains in a fixed complex sector. The local transverse permittivity and the classical Maxwell
maps do not depend on $\diff$.

The main difficulty is to pass from an expansion of the interface
response to an estimate for the coupled scattering problem. Exact
elimination of the longitudinal field gives an interior
Neumann-to-Dirichlet map. Its small-screening-length expansion loses
derivatives and therefore does not give a small perturbation in the
natural trace norm. We overcome this difficulty by combining a uniform
weak estimate for the exact boundary density with the smoothing property
of the far-field map. The same argument gives a second-order expansion
whose coefficients separate curvature and the change in the coupled
Maxwell field. Matching both surface coefficients with a regularized
operator gives the second-order model. The estimates also apply to
observations at a fixed positive distance from the particles and to
plane-wave incidence.

Rigorous effective interface laws are available in other
electromagnetic settings. For a periodic layer of plasmonic
nanoparticles above a conducting plane, Ammari, Li, and
Zou~\cite{AmmariLiZouMetasurfaces} use homogenization to obtain a
Leontovich impedance condition that captures the collective resonances.
Generalized impedance boundary conditions play the same role for
strongly absorbing obstacles \cite{HaddarJolyNguyen}. In both cases, the
condition replaces a whole region. Here, the eliminated scale is the
screening length inside fixed particles, and local Maxwell fields are
retained on both sides of each interface. This also avoids resolving
the screening layer, which volumetric discretizations of the coupled
system must do \cite{MaZhangZou,ChaumontFrelet}.

In nonlocal plasmonics, effective surface descriptions go back to
Feibelman's parameters, which encode moments of induced charge and
current at a metal surface \cite{Feibelman}. Their use in resonance corrections and
generalized Maxwell transmission conditions is developed in
\cite{Christensen,Yang}, and wave-number-dependent parameters describe
stronger tangential confinement \cite{Babaze}. Broadband surface
response measurements are reported in \cite{Chen}. The projected-dipole
model of Yan, Wubs, and Mortensen \cite{YanProjected} represents a
microscopic surface response by a normal dipole layer. The local
analogue model of Luo and collaborators \cite{LuoLAM} instead
represents hydrodynamic screening by a thin dielectric coating on a
local metal, and related planar scattering
formulas are developed in \cite{Pitelet}. Appendix~\ref{app:coating}
gives quantitative estimates for such coatings. Closest to the present
work, Schnitzer et al.\ derive the first-order surface law by formal
matched asymptotics in the quasistatic hydrodynamic model, and identify
a relative error of second order in the screening length
\cite{SchnitzerPRSA}. Our contribution is a rigorous error estimate for
the full retarded scattering problem, including diffusion, together
with a second-order surface approximation.

Recently, Lee, Ruiz, and Yu \cite{LeeRuizYu} developed a spectral theory
for the quasistatic hydrodynamic Drude model on bounded Lipschitz
domains. Their boundary integral formulation gives modal expansions,
describes resonant excitation, and yields perturbation formulas for the
resonances under an assumption of spectral convergence. Here we estimate the retarded far-field operator for
fixed smooth geometries as the screening length tends to zero. For spheres, we prove a pole and
residue comparison for each fixed quasistatic multipole.
For $\diff>0$ the first pole correction contains the diffusion-induced
damping of GNOR.

The paper is organized as follows. Section~\ref{sec:notation} introduces
the geometry, the GNOR model, its dimensionless form, and the asymptotic
regime. Sections~\ref{sec:hdm}--\ref{sec:ntdexp} prove uniform
stability, the exact surface reduction, and the expansion of the surface
operator. Section~\ref{sec:feibelman} relates the reduction to surface
moments. Section~\ref{sec:farfield} proves the far-field estimates, and
Section~\ref{sec:resonance} compares spherical poles and residues.
Appendix~\ref{app:coating} contains the coating estimates and their
numerical verification. Appendix~\ref{app:transmission} gives the
classical transmission estimates used in the proof.

\section{Geometry and governing equations}\label{sec:notation}

\subsection{Geometry and notation}

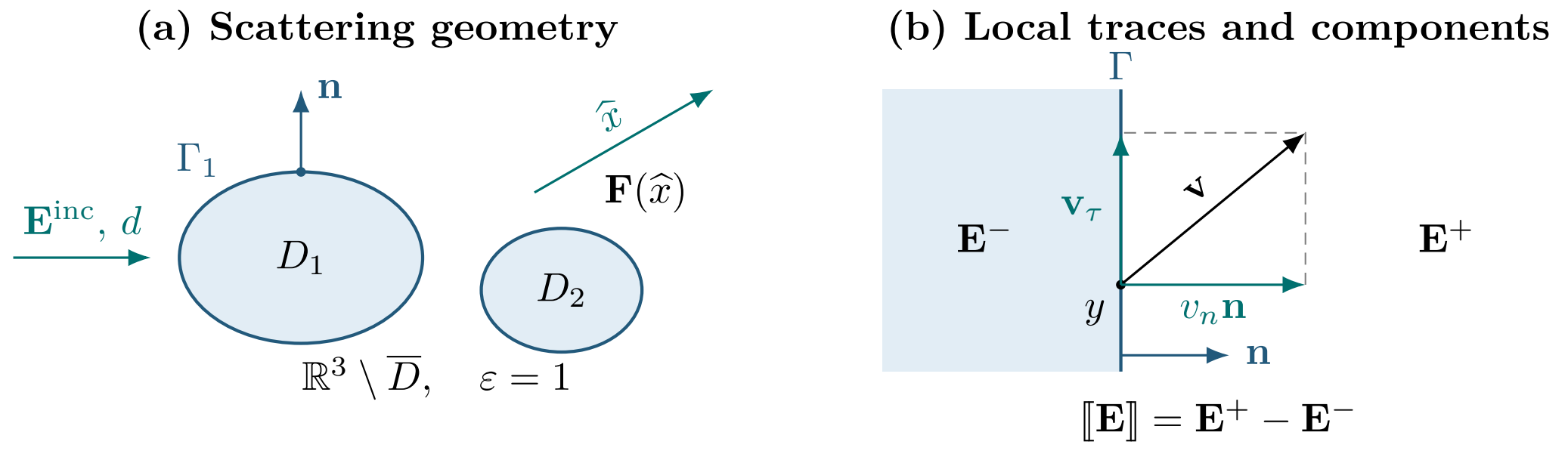
\begin{figure}[!htbp]
\centering
\begin{tikzpicture}[x=0.95cm,y=0.95cm,font=\small]
 \node[font=\small\bfseries] at (3.35,3.45) {(a) Scattering geometry};
 \node[font=\small\bfseries] at (11.1,3.45) {(b) Local traces and components};
 \filldraw[fill=metalwash,figboundary] (2.65,1.35) ellipse (1.12 and 0.79);
 \filldraw[fill=metalwash,figboundary] (5.05,1.05) ellipse (0.74 and 0.57);
 \node at (2.65,1.35) {$D_1$};
 \node at (5.05,1.05) {$D_2$};
 \node at (3.9,0.25) {$\R^3\setminus\overline D,\quad\eps=1$};
 \draw[figarrow,fieldteal] (0,1.35)--(1.27,1.35)
       node[midway,above] {$\E^{\inc},\,d$};
 \fill[metalblue] (2.65,2.14) circle (1.2pt);
 \draw[figarrow,metalblue] (2.65,2.14)--(2.65,2.9)
       node[right] {$\nn$};
 \node[metalblue] at (1.7,2.25) {$\Gamma_1$};
 \draw[figarrow,fieldteal] (4.8,1.95)--(6.45,2.9)
       node[midway,above,sloped] {$\widehat x$};
 \node[anchor=west] at (5.3,1.95) {$\FF(\widehat x)$};
 \fill[metalwash] (8,0.3) rectangle (10.2,2.9);
 \draw[figboundary] (10.2,0.3)--(10.2,2.9);
 \node[metalblue] at (10.2,3.12) {$\Gamma$};
 \node at (8.95,1.55) {$\E^-$};
 \node at (13.2,1.55) {$\E^+$};
 \coordinate (y) at (10.2,1.1);
 \coordinate (v) at (11.9,2.5);
 \fill (y) circle (1.2pt);
 \node[anchor=north east] at (y) {$y$};
 \draw[figarrow] (y)--(v) node[midway,above,sloped] {$\mathbf v$};
 \draw[figarrow,fieldteal] (y)--(11.9,1.1)
       node[midway,below] {$v_n\nn$};
 \draw[figarrow,fieldteal] (y)--(10.2,2.5)
       node[midway,left] {$\mathbf v_\tau$};
 \draw[densely dashed,gray] (11.9,1.1)--(v)--(10.2,2.5);
 \draw[figarrow,metalblue] (10.2,0.45)--(11.2,0.45)
       node[right] {$\nn$};
 \node at (11.1,-0.15) {$\jump{\E}=\E^+-\E^-$};
\end{tikzpicture}
\caption{Cross-sectional sketches of three-dimensional particles.
The normal points into the vacuum; the plus and minus signs indicate
exterior and interior traces. Here $\Gamma_j=\partial D_j$,
$d$ is an incident plane-wave direction, and $\widehat x$ is an observation
direction. The point $y$ lies on $\Gamma$, and the vector decomposition is
$\mathbf v=v_n\nn+\mathbf v_\tau$. In the effective formulation obtained by exact longitudinal reduction,
$\E^-=\ET$ is the divergence-free interior electric field (Proposition~\ref{prop:exact}). All sketches are schematic.}
\label{fig:geometry}
\end{figure}

Let $D=\bigcup_{j=1}^M D_j\subset\R^3$ be a bounded union of disjoint connected
domains with $C^\infty$ boundaries and pairwise disjoint closures.
Put $\Gamma_j=\partial D_j$ and $\Gamma=\partial D$.
The unit normal $\nn$ points from the metal into the vacuum.
Let $\Id$ denote the identity matrix or operator on the space in question.
For a field $\mathbf v$, write
$v_n=\nn\cdot\mathbf v$ and
$\mathbf v_\tau=(\Id-\nn\nn^\top)\mathbf v$.
The jump is $\jump{\mathbf v}=\mathbf v^+-\mathbf v^-$, with $+$ denoting
the vacuum trace. In an effective model the minus trace is that of the
\emph{effective bulk field}; it is not necessarily the microscopic interior
field. We write $\nabla_\Gamma$ and $\diver_\Gamma$ for the surface gradient and
divergence. The Laplace--Beltrami operator
$\Delta_\Gamma=\diver_\Gamma\nabla_\Gamma$ is nonpositive.
The mean curvature $\mathcal H$ is half the sum of the principal curvatures,
and equals $1$ on the unit sphere with this normal convention.
Figure~\ref{fig:geometry} illustrates these conventions.

We write $H^r$ for the usual complex Sobolev spaces; vector-valued
spaces are understood when the unknown is a vector. Spaces on disconnected
boundaries are direct sums. We use $\mathcal L(X,Y)$ for bounded linear
operators from $X$ to $Y$, and a superscript $*$ for the Hilbert adjoint
(or conjugate transpose of a matrix). Norms without a subscript refer to
the spaces specified in the statement. The constants $C,c>0$ may change
from line to line; their dependence on fixed data is indicated in each
result. The notation $A\lesssim B$ means $A\leq CB$, and
$A\asymp B$ means that both $A\lesssim B$ and $B\lesssim A$ hold.
Operator-valued $O(\eta^p)$ remainders are in the stated operator norm.
Inner products
are linear in their first argument; Sobolev duality pairings extend this
convention. The spaces $H(\curl;D)$ and
$H(\diver;D)$ have their usual graph norms, and
\[
 H_0(\diver;D)=\{\mathbf v\in H(\diver;D):\nn\cdot\mathbf v=0
                     \text{ in }H^{-1/2}(\Gamma)\}.
\]

\subsection{The Maxwell--GNOR model}

In the time-harmonic regime, with time dependence $e^{-\ii\omega t}$,
the GNOR model couples the full Maxwell equations to a nonlocal
constitutive equation for the electronic current density $\J$ in the
metal \cite{GNOR}. In the absence of impressed sources near the
particles,
\begin{equation}\label{eq:physmaxwell1}
 \curl\E=\ii\omega\mu_0\HH,
 \qquad
 \curl\HH=\J-\ii\omega\eps_0\E,
\end{equation}
\begin{equation}\label{eq:physcurrent}
 \bigl(\beta^2+D_{\mathrm{diff}}(\gamma-\ii\omega)\bigr)\nabla\diver\J
   +\omega(\omega+\ii\gamma)\J
   =\ii\omega\eps_0\omega_p^2\E\quad\text{in the metal},
\end{equation}
\begin{equation}\label{eq:physbc}
 \nn\times\jump{\E}=0,\qquad
 \nn\times\jump{\HH}=0,\qquad
 \J\cdot\nn=0\quad\text{on the metal boundary}.
\end{equation}
Here $\E$ and $\HH$ are the electric and magnetic fields, $\J$ is
extended by zero outside the metal, $\omega$ is the angular frequency,
and $\ii^2=-1$. The constants $\eps_0$ and $\mu_0$ are the vacuum
permittivity and permeability, and $c=(\eps_0\mu_0)^{-1/2}$ and
$Z_0=(\mu_0/\eps_0)^{1/2}$. The metal is characterized by its plasma
frequency $\omega_p>0$, damping rate $\gamma>0$, hydrodynamic velocity
$\beta>0$, and a spatially constant, frequency-independent diffusion
coefficient $D_{\mathrm{diff}}\geq0$. The host and bound-electron background
relative permittivities are $1$, and the relative permeability is $1$
everywhere.

For $D_{\mathrm{diff}}=0$, \eqref{eq:physcurrent} is the hydrodynamic Drude
current equation. The GNOR model adds diffusion by replacing $\beta^2$
with $\beta^2+D_{\mathrm{diff}}(\gamma-\ii\omega)$; see also
\cite[Eq.~(17)]{StamatopoulouTserkezis}. The last condition in
\eqref{eq:physbc} is a hard wall. Electrons do not cross the metal
boundary, so spill-out and interparticle charge transfer are excluded.
Here $\J$ is the \emph{total} electronic current, including diffusion,
and the hard-wall condition is imposed on this total current.

\subsection{Dimensionless formulation}

Let $a>0$ be a reference particle length, and let the metal occupy
$aD$, with $D$ as above. Fix an electric-field
scale $E_*>0$ and set
\begin{equation}\label{eq:scaling}
 \begin{gathered}
 x'=\frac xa,\qquad z=\frac\omega{\omega_p},\qquad
 \gamma'=\frac\gamma{\omega_p},\qquad
 \eta=\frac{\beta}{a\omega_p},\qquad
 \sigma=\frac{a\omega_p}{c},\qquad \kappa=\sigma z,\\
 \E'=\frac{\E}{E_*},\qquad \HH'=\frac{Z_0\HH}{E_*},\qquad
 \J'=\frac{\J}{\eps_0\omega_pE_*},
 \end{gathered}
\end{equation}
together with the dimensionless diffusion parameter
\begin{equation}\label{eq:gnorparameters}
 \diff=\frac{\omega_pD_{\mathrm{diff}}}{\beta^2},
 \qquad\text{so that}\qquad
 \frac{D_{\mathrm{diff}}}{a^2\omega_p}=\eta^2\diff.
\end{equation}
From now on we drop the primes. By this abuse of notation, $x$,
$\gamma$, $\E$, $\HH$, and $\J$ denote dimensionless quantities, and
all lengths, surface differential operators, and surface-response
parameters are dimensionless unless otherwise indicated.
Thus $\ell=\beta/\omega_p$ is the reference screening length and
$\eta=\ell/a$. The hydrodynamic Drude model is the case $\diff=0$, and
$\eta$ remains a real positive reference length ratio for every $\diff$.
In these variables, with $\J$ extended by zero outside $D$,
\eqref{eq:physmaxwell1}--\eqref{eq:physbc} become
\begin{equation}\label{eq:maxwell1}
 \curl\E=\ii\kappa\HH,\qquad
 \curl\HH=\sigma\J-\ii\kappa\E\quad\text{in }\R^3,
\end{equation}
\begin{equation}\label{eq:hdcurrent}
 \eta^2\bigl(1+\diff(\gamma-\ii z)\bigr)\nabla\diver\J
   +z(z+\ii\gamma)\J=\ii z\E\quad\text{in }D,
\end{equation}
and
\begin{equation}\label{eq:physicalbc}
 \nn\times\jump{\E}=0,\qquad
 \nn\times\jump{\HH}=0,\qquad
 \J\cdot\nn=0\quad\text{on }\Gamma.
\end{equation}
For $z\ne0$, the divergence relations follow from the curl
equations:
\[
 \diver\HH=0,\qquad
 \diver\E=\frac{1}{\ii z}\diver\J\ \text{in }D,
 \qquad \diver\E=0\ \text{outside }\overline D.
\]
Together with \eqref{eq:physicalbc}, these equations also imply
$\jump{E_n}=\jump{H_n}=0$ for the physical fields, in the trace sense.

The incident electric field is a Herglotz wave
\begin{equation}\label{eq:herglotz}
 \E^{\inc}(x)=\int_\Sph \mathbf a(d)e^{\ii\kappa d\cdot x}\dd S(d),
 \qquad \mathbf a\in L_t^2(\Sph)
 :=\{\mathbf a\in L^2(\Sph)^3:d\cdot\mathbf a(d)=0\}.
\end{equation}
The corresponding magnetic field has density $d\times\mathbf a(d)$.
In the exterior, write $\E=\E^{\inc}+\E^{\scat}$ and
$\HH=\HH^{\inc}+\HH^{\scat}$. The scattered fields satisfy the
Silver--M\"uller radiation condition
\begin{equation}\label{eq:silvermuller}
 \lim_{r\to\infty}r\bigl(
 \HH^{\scat}(r\widehat x)\times\widehat x
       -\E^{\scat}(r\widehat x)\bigr)=0,
 \qquad \widehat x\in\Sph,
\end{equation}
uniformly in $\widehat x$. The outgoing electric far field is defined by
\begin{equation}\label{eq:farfielddef}
 \E^{\scat}(r\widehat x)=\frac{e^{\ii\kappa r}}r
       \FF(\widehat x)+O(r^{-2}),\qquad r\longrightarrow\infty.
\end{equation}
In physical units, the electric far field is $aE_*\FF$.

An outgoing weak GNOR solution consists of
$\E,\HH\in H_{\mathrm{loc}}(\curl;\R^3)$ and
$\J\in H_0(\diver;D)$ satisfying
\eqref{eq:maxwell1}--\eqref{eq:hdcurrent} in the distributional sense,
\eqref{eq:physicalbc} in the trace sense, and
\eqref{eq:silvermuller} for the scattered fields. Here
$H_{\mathrm{loc}}(\curl;\R^3)$ means $H(\curl)$ regularity on every
bounded ball.

\subsection{Material functions and the asymptotic regime}

With square roots of positive real part, define
\begin{equation}\label{eq:material}
 \Omega(z)=z(z+\ii\gamma),\qquad
 \eps(z)=1-\Omega(z)^{-1},\qquad
 s(z)=\sqrt{1-\Omega(z)},
\end{equation}
and
\begin{equation}\label{eq:screening}
 \varg(z)=1+\diff(\gamma-\ii z),\qquad
 \sg(z)=\sqrt{\frac{1-\Omega(z)}{\varg(z)}},
 \qquad \Qg(z)=\frac{\sg(z)}\eta.
\end{equation}
Thus, $s_0=s$, and \eqref{eq:hdcurrent} reads
$\eta^2\varg\nabla\diver\J+\Omega\J=\ii z\E$.
The two square roots play different roles. The transverse identity
$1-\eps^{-1}=s^{-2}$ involves $s$, whereas $\sg$ governs the
longitudinal screening. They coincide only at $\diff=0$.
Whenever $\eps$ occurs inside a jump, its exterior value is understood
to be $1$.

Throughout, the domains and their positive separations are fixed,
and
\begin{equation}\label{eq:regime}
 \gamma>0\text{ fixed},\quad I=[z_-,z_+]\Subset(0,1),\quad
 \sigma>0\text{ fixed},\quad 0\leq\diff\leq\diff_*<\infty,
 \quad \eta\downarrow0.
\end{equation}
On $I$, the quantities $s$, $s^{-1}$, and $\eps^{-1}$ are bounded.
The screening parameter obeys the sector identity
\begin{equation}\label{eq:gnorsector}
 \operatorname{Re}\sg(z)^2
 =\frac{1-z^2+\diff\gamma}
        {(1+\diff\gamma)^2+\diff^2z^2}>0.
\end{equation}
Compactness of $I\times[0,\diff_*]$ therefore gives $|\sg|\asymp1$ and
$\operatorname{Re}\Qg^2\geq c|\Qg|^2>0$, with constants independent of
$\eta$. All constants below may depend on $\diff_*$ but are uniform in
$z\in I$ and $\diff\in[0,\diff_*]$. We suppress the dependence on
$\diff$ in the notation unless it is needed.
Since $\eta\sigma=\beta/c$, this is a small-screening-length limit at
fixed electromagnetic size, rather than a change of particle size with
all material parameters fixed. By \eqref{eq:gnorparameters}, it makes
both the pressure and the diffusion lengths small. At fixed
$a,\omega_p$ and fixed positive $D_{\mathrm{diff}}$, letting
$\beta\to0$ violates \eqref{eq:regime} and does not give the local
Drude limit.

\begin{remark}\label{rem:unboundeddiff}
The bound $\diff_*$ controls only the size of $\sg$. The ratio
$\operatorname{Re}\sg^2/|\sg^2|$ has the positive limit
$\gamma/(|1-\Omega|\sqrt{\gamma^2+z^2})$ as $\diff\to\infty$, so the
sector is uniform for all $\diff\geq0$ when $\gamma>0$. In contrast,
$|\sg|\asymp(1+\diff)^{-1/2}$. For unbounded $\diff$ the natural small
parameter is $\eta/|\sg|\asymp\eta(1+\diff)^{1/2}$, where
$\eta(1+\diff)^{1/2}$ combines the pressure and diffusion length ratios
$\eta$ and $\eta\sqrt{\diff}$ in quadrature. For the regularization in
Section~\ref{sec:farfield}, the contraction estimate there is
$O(\sqrt{\eta}/|\sg|)$. Smallness of $\eta/|\sg|$ alone does not ensure
that this bound tends to zero.
\end{remark}

\section{Uniform stability of the GNOR problem}\label{sec:hdm}

In this section, we use the volume representation to obtain a stability
estimate independent of the screening length. Let
$\Phi_\kappa(x)=e^{\ii\kappa|x|}/(4\pi|x|)$ and let $\GG$ be the
outgoing electric volume potential, restricted to $D$ when appropriate:
\begin{equation}\label{eq:green}
 \GG\mathbf p=\left(\Id+\kappa^{-2}\nabla\nabla\cdot\right)
                      (\Phi_\kappa*\mathbf p).
\end{equation}
Sources supported in $D$ are extended by zero. The derivatives are distributional. This is a bounded operator on
$L^2(D)^3$ for fixed bounded $D$ and $\kappa\neq0$ in the frequency sets
used below. Here $*$ denotes convolution, whereas an operator superscript $*$ denotes
the adjoint. For an integrable dimensionless polarization density $\mathbf p$
of compact support, the field $\kappa^2\GG\mathbf p$ has the far-field
operator
\begin{equation}\label{eq:radiation}
 (\mathcal C_z\mathbf p)(\widehat x)
 =\frac{\kappa^2}{4\pi}(\Id-\widehat x\widehat x^\top)
       \int_{\R^3} e^{-\ii\kappa\widehat x\cdot y}\mathbf p(y)\dd y.
\end{equation}

We establish an $\eta$-uniform estimate for the scattering problem
\eqref{eq:maxwell1}--\eqref{eq:hdcurrent}, with
\eqref{eq:physicalbc} and \eqref{eq:silvermuller}, in the regime
\eqref{eq:regime}.

\begin{lemma}[Uniform dissipative estimate]\label{lem:energy}
For each $z\in I$, $\eta>0$, and $\diff\geq0$,
\eqref{eq:maxwell1}--\eqref{eq:hdcurrent} has a unique
outgoing weak solution. Uniformly for $z\in I$ and $\diff\in[0,\diff_*]$,
it satisfies
\begin{equation}\label{eq:currentbound}
 \norm{\J}_{L^2(D)}\leq\gamma^{-1}\norm{\E^{\inc}}_{L^2(D)},
 \qquad \norm{\E}_{L^2(D)}\leq C\norm{\mathbf a}_{L^2(\Sph)}.
\end{equation}
\end{lemma}
\begin{proof}
Extend $\J$ by zero. The electric field is
\begin{equation}\label{eq:currentgreen}
 \E=\E^{\inc}+\ii\sigma^2 z\GG\J.
\end{equation}
On $V=H_0(\diver;D)$ the current equation is $a_\eta(\J,v)
=(\E^{\inc},v)$, where
\[
 a_\eta(u,v)=(\gamma-\ii z)(u,v)
 +\eta^2\left(\diff+\frac{\ii(1+\diff\gamma)}z\right)
             (\diver u,\diver v)
 -\ii\sigma^2z(\GG u,v).
\]
For real $\kappa>0$, the outgoing radiation identity is
\[
 \operatorname{Im}(\GG u,u)_{L^2(D)}
 =\kappa^{-3}\|\mathcal C_z u\|_{L^2(\Sph)}^2\geq0.
\]
It follows either from Green's formula or from the Fourier representation
of $\operatorname{Im}\Phi_\kappa$. Thus
\begin{equation}\label{eq:gnordissipation}
 \operatorname{Re}a_\eta(u,u)
 =\gamma\|u\|^2+\eta^2\diff\|\diver u\|^2
       +\sigma^2z\operatorname{Im}(\GG u,u)
 \geq\gamma\|u\|^2+\eta^2\diff\|\diver u\|^2.
\end{equation}
The imaginary coefficient of the divergence term is
$\eta^2(1+\diff\gamma)/z\geq\eta^2/z$. The other terms in the imaginary
part are bounded in absolute value by $C_I\norm u^2$.
For a sufficiently small fixed $\theta>0$, chosen independently of
$\diff\geq0$,
\[
 \operatorname{Re}\{e^{-\ii\theta}a_\eta(u,u)\}
 \geq c\bigl(\norm u^2+\eta^2\norm{\diver u}^2\bigr).
\]
For fixed $\eta>0$, the weighted norm on the right is equivalent to
the graph norm of $V$. Lax--Milgram gives existence and uniqueness in
$V$, and
\eqref{eq:currentgreen} reconstructs the outgoing Maxwell fields.
Taking real parts with $v=\J$ in \eqref{eq:gnordissipation} gives the
first bound in \eqref{eq:currentbound}. Boundedness of $\GG$ gives the
second.
\end{proof}

Diffusion only adds the nonnegative term $\eta^2\diff\|\diver u\|^2$
to the dissipation. Uniqueness for hydrodynamic media with complex
constitutive parameters, including GNOR, is also studied in
\cite{Mystilidis}. The additional point here is the estimate uniform in
the screening length and the diffusion parameter.

\section{Exact elimination of the longitudinal field}\label{sec:exact}

We separate the divergence-free electromagnetic field from the
longitudinal screening field and eliminate the latter exactly.
Let $\mathcal N_Q:H^{-1/2}(\Gamma)\to H^{1/2}(\Gamma)$ be the interior
Neumann-to-Dirichlet map for the modified Helmholtz equation:
\begin{equation}\label{eq:ntd}
 (\Delta-Q^2)\psi=0\text{ in }D,\qquad
 \partial_n\psi=f\text{ on }\Gamma,\qquad
 \mathcal N_Qf=\psi|_\Gamma.
\end{equation}
The map is defined componentwise on $D_j$. Set $\NN=\mathcal N_{\Qg}$,
with $\Qg=\sg(z)/\eta$ from \eqref{eq:screening}.

For the reduced surface problems in this section and in
Section~\ref{sec:farfield}, an energy solution
consists of local Maxwell fields in $H(\curl;D)$ and in
$H(\curl;B_R\setminus\overline D)$ for every ball $B_R$ containing
$\overline D$, together with a scalar interface potential
$v\in H^{1/2}(\Gamma)$ satisfying the stated surface equations and
transmission conditions. The exterior scattered fields satisfy
\eqref{eq:silvermuller}. In particular, the effective electric field
has piecewise $H(\curl)$ regularity; its prescribed tangential jump
need not vanish.

\begin{proposition}[Exact surface reduction]\label{prop:exact}
In the regime \eqref{eq:regime}, define in $D$
\begin{equation}\label{eq:split}
 \psi=\Qg^{-2}\diver\E,\qquad \ET=\E-\nabla\psi.
\end{equation}
Then the exterior field agrees exactly with the solution of the local bulk
Maxwell equations
\begin{equation}\label{eq:transverse}
 \curl\ET=\ii\kappa\HH,\qquad
 \curl\HH=-\ii\kappa\eps\ET\quad\text{in }D,
\end{equation}
coupled to exterior vacuum Maxwell equations by
\begin{align}
 \E_\tau^+-\ET_\tau&=\nabla_\Gamma v,
 &v&=\NN f, & f&=(\eps-1)E_n^{\mathrm{T}},\label{eq:exactbc}\\
 \HH_\tau^+&=\HH_\tau^-,
 &E_n^+&=\eps E_n^{\mathrm{T}},
 &H_n^+&=H_n^-.\nonumber
\end{align}
Conversely, a solution of this effective problem reconstructs the
GNOR solution by solving \eqref{eq:ntd} and setting
\begin{equation}\label{eq:reconstruct}
 \E^-_{\gnor}=\ET+\nabla\psi,\qquad
 \J=\ii z\bigl((1-\eps)\ET+\nabla\psi\bigr).
\end{equation}
The form of \eqref{eq:transverse} and \eqref{eq:reconstruct} is
independent of $\diff$. The potential $\psi$ depends on $\diff$ through $\Qg$.
\end{proposition}
\begin{proof}
Taking divergence in \eqref{eq:maxwell1} gives
$\diver\J=\ii z\diver\E$. For each fixed $\eta>0$,
\eqref{eq:hdcurrent} also gives
\[
 \nabla\diver\J=\eta^{-2}\varg^{-1}(\ii z\E-\Omega\J)\in L^2(D)^3,
\]
since $\varg\neq0$.
Since $\J\in H_0(\diver;D)$, this implies
$\diver\E\in H^1(D)$ and justifies $\nabla\psi\in L^2(D)^3$.
Taking the divergence of \eqref{eq:hdcurrent} now yields
$[\eta^2\varg\Delta+\Omega-1]\diver\E=0$.
It follows that $(\Delta-\Qg^2)\psi=0$ and $\diver\ET=0$.
Substitution in the current equation gives \eqref{eq:reconstruct}, since
$\eta^2\varg\nabla\diver\J=\ii z(1-\Omega)\nabla\psi$ and
$\Omega^{-1}=1-\eps$. This gives \eqref{eq:transverse}. The hard-wall condition becomes
$\partial_n\psi=(\eps-1)E_n^{\mathrm{T}}$.
Tangential continuity of the physical electric field gives
\eqref{eq:exactbc}. Taking normal components of the curl equations gives
the two normal transmission conditions. All steps can be reversed.
\end{proof}

\begin{corollary}\label{cor:weaktrace}
The exact boundary density $f_\eta=(\eps-1)E_n^{\mathrm{T}}$ satisfies
\begin{equation}\label{eq:weaktrace}
 \sup_{z\in I,\ 0\leq\diff\leq\diff_*}\norm{f_\eta}_{H^{-1/2}(\Gamma)}
       \leq C\norm{\mathbf a}_{L^2(\Sph)},
\end{equation}
with $C$ independent of $\eta$.
\end{corollary}
\begin{proof}
By \eqref{eq:reconstruct},
$\ET=\eps^{-1}(\E-\J/(\ii z))$ is uniformly bounded in $L^2(D)^3$.
It is divergence free. In particular, this argument does not require
a uniform bound on $\diver\E$ or on $\nabla\psi$ separately.
The normal trace theorem in $H(\diver;D)$ proves
\eqref{eq:weaktrace}.
\end{proof}

The reduction is exact before taking $\eta\to0$. For disjoint particles
the electronic map $\NN$ is block diagonal, whereas the exterior Maxwell
problem couples all components. Analytic continuation preserves this
equivalence only in frequency regions where the eliminated scalar
Neumann problem is invertible and the factors in \eqref{eq:split} are
defined.

\section{Expansion of the surface operator}\label{sec:ntdexp}

We next expand the Neumann-to-Dirichlet map on smooth data. The derivative
loss in the remainder will be retained in the exterior error estimate.

\begin{lemma}[Screening-layer expansion]\label{lem:ntdexp}
Let $\Gamma$ be fixed and smooth. In a closed sector
$\operatorname{Re}Q>0$, $\operatorname{Re}Q^2\geq c|Q|^2>0$,
for $|Q|$ sufficiently large and
each integer $r\geq0$,
\begin{equation}\label{eq:ntdexp}
 \norm{\mathcal N_Qf-Q^{-1}f-Q^{-2}\mathcal H f}_{H^r(\Gamma)}
 \leq C_r|Q|^{-3}\norm f_{H^{r+3}(\Gamma)}.
\end{equation}
In particular,
\begin{align}
 \norm{\NN f}_{H^r}&\leq C_r\eta\norm f_{H^{r+3}},\label{eq:Nsmall}\\
 \norm{(\NN-\eta\sg^{-1}\Id)f}_{H^r}
     &\leq C_r\eta^2\norm f_{H^{r+3}},\label{eq:Napprox}\\
 \left\|\left(\NN-\frac\eta{\sg}\Id
       -\frac{\eta^2}{\sg^2}\mathcal H\right)f\right\|_{H^r}
     &\leq C_r\eta^3\norm f_{H^{r+3}}.\label{eq:gnorntdexp}
\end{align}
The constants depend on $\Gamma$, $r$, and the fixed sector. By
\eqref{eq:gnorsector}, the estimates for $\NN$ are uniform in
$z\in I$ and $\diff\in[0,\diff_*]$.
The same estimates hold for the $L^2(\Gamma)$ adjoints.
\end{lemma}
\begin{proof}
In an interior collar, write $x=y-t\nn(y)$, $t\geq0$, and choose a smooth
cutoff $\chi(t)$ equal to one near zero. The outward derivative is
$-\partial_t$ at $t=0$, and
\[
 \Delta=\partial_t^2-2\mathcal H\partial_t+\Delta_\Gamma
                 +O(t)\partial_t+O(t)\partial_y^2+O(t)\partial_y.
\]
Here, the $O(t)$ symbols denote smooth differential-operator coefficients
vanishing at $t=0$. Define
\[
 \psi_{\mathrm{app}}(y,t)=\chi(t)e^{-Qt}
       \left(\frac{f(y)}Q+\frac{\mathcal H(y)f(y)}{Q^2}
                         +\frac{t\mathcal H(y)f(y)}Q\right).
\]
Its normal derivative at zero is exactly $f$.
The leading curvature residual cancels, and direct differentiation gives
\[
 \norm{(\Delta-Q^2)\psi_{\mathrm{app}}}_{L^2(D)}
       \leq C|Q|^{-3/2}\norm f_{H^2(\Gamma)}.
\]
The terms involving the cutoff are exponentially small. Set
$R=(\Delta-Q^2)\psi_{\mathrm{app}}$ and $w=\psi-\psi_{\mathrm{app}}$.
Then $(\Delta-Q^2)w=-R$ and the Neumann data of $w$ vanish. The parameter-dependent
energy estimate and the multiplicative trace inequality give
\begin{align*}
 \norm w_{L^2(D)}&\leq C|Q|^{-2}\norm R_{L^2(D)},
 &\norm{\nabla w}_{L^2(D)}&\leq C|Q|^{-1}\norm R_{L^2(D)},\\
 \norm{w|_\Gamma}_{L^2(\Gamma)}&\leq C|Q|^{-3/2}\norm R_{L^2(D)}.
\end{align*}
This proves the $r=0$ estimate, even with $H^2$ data.

For completeness, use tangential derivatives in the collar to avoid
introducing powers of $Q$ by differentiating $e^{-Qt}$ normally.
Let $\|\cdot\|_{r,\tau}$ be the sum of the $L^2$ norms of derivatives
of order at most $r$ in the boundary variables in a fixed partition of
unity; away from the collar, use ordinary interior charts. Direct
differentiation of the same residual gives
\[
 \|R\|_{r,\tau}\leq C_r|Q|^{-3/2}\|f\|_{H^{r+2}(\Gamma)}.
\]
The Neumann resolvent satisfies the tangential parameter estimates
\begin{equation}\label{eq:tangentialresolvent}
 |Q|^2\|w\|_{r,\tau}+|Q|\|\nabla w\|_{r,\tau}
 \leq C_r\|R\|_{r,\tau}.
\end{equation}
To obtain these estimates, transport the boundary coordinate vector
fields along the normal lines, so that they commute with $\partial_t$
and preserve the homogeneous Neumann condition. Apply tangential
difference quotients to the variational equation. The commutators of
the divergence-form principal part are bounded by the weighted energy
norms of lower tangential derivatives. The zero-order parameter $Q^2$
commutes with differentiation. Testing with the differentiated solution
and using $\operatorname{Re}Q^2\geq c|Q|^2$ gives, if
$E_r^2=|Q|^2\|w\|_{r,\tau}^2+\|\nabla w\|_{r,\tau}^2$,
\[
 E_r\leq C_r\bigl(|Q|^{-1}\|R\|_{r,\tau}+E_{r-1}\bigr),
 \qquad E_0\leq C|Q|^{-1}\|R\|_{0,\tau}.
\]
Cutoff commutators have the same lower-order bound. Induction proves
\eqref{eq:tangentialresolvent}. This is the usual tangential proof of
parameter elliptic regularity for the smooth Neumann problem.

The multiplicative trace inequality applied to these derivatives yields
\[
 \|w|_\Gamma\|_{H^r}
 \leq C_r\|w\|_{r,\tau}^{1/2}
       (\|\nabla w\|_{r,\tau}+\|w\|_{r,\tau})^{1/2}
 \leq C_r|Q|^{-3}\|f\|_{H^{r+2}}.
\]
The weaker bound with $H^{r+3}$ data in \eqref{eq:ntdexp} follows.
Density extends the estimate from smooth data. Green's identity gives
$\mathcal N_Q^*=\mathcal N_{\overline Q}$; the conjugate parameter lies
in the same sector, proving the adjoint estimates.
\end{proof}

On the plane, let $\xi\in\R^2$ denote the tangential Fourier variable
and let a hat denote the tangential Fourier transform. The exact symbol is
\begin{equation}\label{eq:planarN}
 \widehat{\mathcal N_Qf}(\xi)
       =(Q^2+|\xi|^2)^{-1/2}\widehat f(\xi).
\end{equation}
Thus the expansion requires slow tangential variation relative to
$|Q|^{-1}$. In particular, \eqref{eq:Napprox} is not an operator-norm
estimate $H^{-1/2}\to H^{1/2}$ with error $O(\eta^2)$.
Derivative loss cannot be suppressed when the boundary law is solved
self-consistently.

\section{Surface moments and the Feibelman law}\label{sec:feibelman}

We relate the surface reduction to the moment description of induced
charge and current.

At a locally planar interface, let $t$ be signed distance, positive into
the vacuum, and let $\rho(t)$ be the
induced dimensionless charge density. Assume that the following moment
integrals converge and that the current excess decays sufficiently fast
for the stated integration by parts. For a chosen excitation define
\[
 q_s=\int_\R\rho(t)\dd t,\qquad
 \pi_s=\int_\R t\rho(t)\dd t,\qquad
 d_\perp=\pi_s/q_s,
\]
whenever $q_s\neq0$. If $\mathbf j_\tau(t)$ approaches constant bulk
values on either side, let $\mathbf j_\tau^{\mathrm{step}}$ be the
piecewise constant function equal to those limits for $t<0$ and $t>0$,
respectively. Define the tangential excess current
\[
 \mathbf k_s=\int_\R(\mathbf j_\tau-\mathbf j_\tau^{\mathrm{step}})\dd t.
\]
Integration by parts gives
$\mathbf k_s=-\int t\partial_t\mathbf j_\tau\dd t$.
For an in-plane isotropic, reflection-invariant response this is written
$\mathbf k_s=-d_\parallel\jump{\mathbf j_\tau}$.
At leading order $\jump{\E_\tau}=0$, so
$\jump{\mathbf j_\tau}=-\ii z\jump{\eps\E_\tau}$.
In this section, jumps refer to extrapolated effective bulk fields.
In particular, they are not the jumps of the resolved microscopic fields
in \eqref{eq:physicalbc}. The identification $q_s=\jump{E_n}$ is used
at leading order in a locally planar, slowly varying cell. Consequently,
to the order retained in the surface law,
\begin{equation}\label{eq:surfaceconstitutive}
 \pi_s=d_\perp\jump{E_n},\qquad
 \mathbf k_s=\ii z d_\parallel\jump{\eps\E_\tau}.
\end{equation}
The resulting first-order surface conditions are
\begin{align}
 \jump{\E_\tau}&=-\nabla_\Gamma(d_\perp\jump{E_n}),\label{eq:feibelman}\\
 \nn\times\jump{\HH}&=\ii\kappa d_\parallel\jump{\eps\E_\tau},\nonumber\\
 \jump{\eps E_n}&=\diver_\Gamma(d_\parallel\jump{\eps\E_\tau}),
 &\jump{H_n}&=0.\nonumber
\end{align}
They follow from the electric dipole-layer jump, Ampere's law with surface
current $\sigma\mathbf k_s$, and surface charge conservation. Products
must remain inside the surface derivatives if the coefficients vary.
These conventions agree with the surface-response formulation in
\cite{Feibelman,Yang}.

\begin{figure}[htbp]
\centering
\begin{tikzpicture}[x=0.95cm,y=0.95cm,font=\small]
 \node[font=\small\bfseries] at (3.05,3.55) {(a) Resolved screening layer};
 \node[font=\small\bfseries] at (11.2,3.55) {(b) Effective interface};
 \fill[metalwash] (0.4,0.65) rectangle (4.3,2.85);
 \fill[screentint] (3.3,0.65) rectangle (4.3,2.85);
 \draw[figboundary] (4.3,0.65)--(4.3,2.9);
 \draw[figarrow] (0.35,0.65)--(6,0.65) node[right] {$\Qg t$};
 \draw[screenorange,very thick,domain=-3.8:0,samples=80]
       plot ({4.3+\x},{0.65+1.8*exp(\x)});
 \draw[screenorange,thick] (4.3,0.65)--(5.7,0.65);
 \draw[screenorange,densely dashed] (3.3,0.65)--(3.3,2.65);
 \foreach \q in {-3,-2,-1,0}
   {\draw ({4.3+\q},0.6)--({4.3+\q},0.7);
    \node[below] at ({4.3+\q},0.57) {$\q$};}
 \node[anchor=east] at (3.9,3.12) {$\rho/\rho(0^-)$};
 \node at (1.55,2.55) {metal};
 \node at (5.15,2.55) {vacuum};
 \node[screenorange] at (2.9,-0.1) {centroid: $t=d_\perp=-\Qg^{-1}$};
 \draw[figarrow] (6.45,1.9)--(7.55,1.9);
 \node[align=center,font=\scriptsize] at (7,2.5) {eliminate\\$\psi$};
 \fill[metalwash] (8,0.65) rectangle (10.4,2.85);
 \draw[screenorange,very thick] (10.4,0.65)--(10.4,2.85);
 \node at (9.1,1.9) {$\ET$};
 \node at (12.5,1.9) {$\E^+$};
 \node[screenorange,anchor=west] at (10.6,2.8) {$v=\psi|_\Gamma$};
 \draw[figarrow,metalblue] (10.4,1)--(11.4,1) node[right] {$\nn$};
 \node at (11.2,0.15) {$v=\mathcal N_{\Qg}f,\quad f=(\eps-1)E_n^{\mathrm{T}}$};
 \node at (11.2,-0.4) {$\jump{\E_\tau}=\nabla_\Gamma v$};
\end{tikzpicture}
\caption{Screening and surface response. The left panel shows the planar,
zero tangential wave-number profile $\rho(t)/\rho(0^-)=e^{\Qg t}$ for $t<0$,
with zero charge in the vacuum. The plot is drawn for $\diff=0$ at a
lossless subplasma frequency, where $\Qg>0$. For $\diff>0$ or
$\gamma>0$ the decay rate $\Qg$ is generally complex.
Here $t$ is signed outward distance, opposite to the nonnegative inward
coordinate used in Lemma~\ref{lem:ntdexp}. The exact surface map on the
right retains the longitudinal field through $v=\mathcal N_{\Qg}f$.
Its leading approximation $v=\Qg^{-1}f$ gives $d_\perp=-\Qg^{-1}$ and
$d_\parallel=0$.}
\label{fig:screening}
\end{figure}
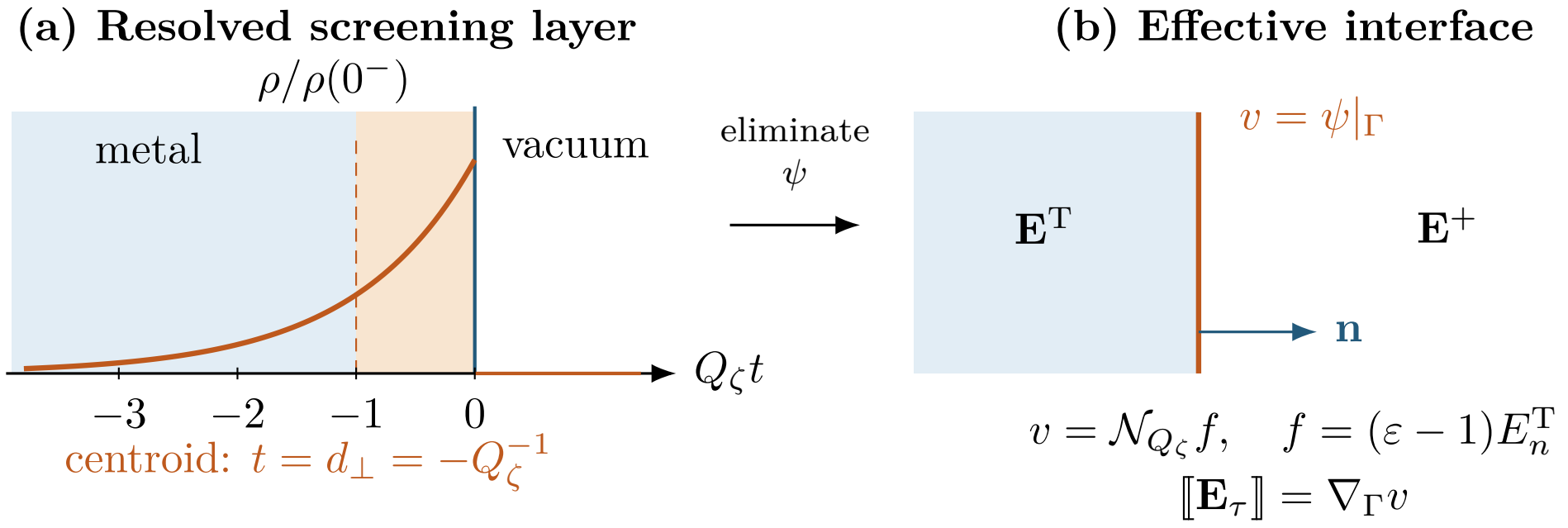

For the hard-wall GNOR model, \eqref{eq:exactbc} gives
$\jump{E_n}=(\eps-1)E_n^{\mathrm{T}}=f$; hence
\begin{equation}\label{eq:hdmd}
 d_\perp=-\eta/\sg,\qquad d_\parallel=0
\end{equation}
to first order. On a plane the full normal response is instead the
operator with symbol
$d_\perp(z,\xi)=-(\Qg^2+|\xi|^2)^{-1/2}$.
In particular, a local surface coefficient is already an asymptotic
replacement of a spatially dispersive response; compare~\cite{Babaze}.
Through $\sg$, the coefficient $d_\perp$ also encodes the diffusion
contribution. The link between GNOR and microscopic surface-response
parameters is discussed in \cite{SvendsenGNOR}. No such microscopic
derivation is assumed below.
Figure~\ref{fig:screening} illustrates the screening profile and its
elimination in favor of a surface unknown.

\section{A controlled far-field approximation}\label{sec:farfield}

We use a regularization to give a well-posed first-order surface model on
the full trace space. Define
\begin{equation}\label{eq:regN}
 R_\eta=(\Id-\eta\Delta_\Gamma)^{-1},\qquad
 \mathcal N_\eta^{\eff}=\frac\eta{\sg}R_\eta.
\end{equation}
The effective fields solve the bulk equations \eqref{eq:transverse} and
exterior Maxwell equations with the interface conditions
\begin{equation}\label{eq:regsurface}
 \jump{\E_\tau}=\nabla_\Gamma v,\qquad
 (\Id-\eta\Delta_\Gamma)v=\frac\eta{\sg}(\eps-1)E_n^-,\qquad
 \jump{\HH_\tau}=0,\quad\jump{\eps E_n}=0,\quad\jump{H_n}=0.
\end{equation}
All interior traces in
\eqref{eq:regsurface} refer to the effective local bulk field.
The regularizer $R_\eta$ uses the real parameter $\eta$ for every
$\diff$, and diffusion enters only through the scalar $\sg$.
The regularization is a mathematical choice. It changes the surface law
only at second order on smooth particle-scale data, but is not asserted
to be a microscopic high-wave-number law or a passive model at all
frequencies.

For numerical implementation, $v$ can be retained as a scalar unknown
on the surface mesh, coupled to a local Maxwell discretization through
\eqref{eq:regsurface}. Thus the additional electronic degrees of freedom
are confined to $\Gamma$, and no volumetric longitudinal field needs to
be discretized. The scalar surface equation is coupled to the Maxwell
field through $E_n^-$; it is not an independent postprocessing step.
Evaluating the exact map $\NN$ in \eqref{eq:exactbc} may still require
costly interior solves. Likewise, explicitly meshing the thin coatings
of Appendix~\ref{app:coating} can reintroduce fine normal resolution.
The computational advantage is therefore clearest for the direct
surface formulation; its magnitude depends on the discretization and
solver, including the volumetric method used for comparison.

We record the classical transmission maps used in the proof.
The reference problem is the local Drude transmission problem. Its
fields $\E_0$ and $\HH_0$ satisfy
\begin{equation}\label{eq:localmaxwell}
 \curl\E_0=\ii\kappa\HH_0,\qquad
 \curl\HH_0=-\ii\kappa\eps_D\E_0,
 \qquad
 \eps_D(x,z)=
 \begin{cases}\eps(z),&x\in D,\\1,&x\notin\overline D,\end{cases}
\end{equation}
with continuous tangential components, the incident field
\eqref{eq:herglotz}, and the radiation condition
\eqref{eq:silvermuller}. No hard-wall condition is imposed, since the
limit $\eta\downarrow0$ is singular at the interface. Let
$f_0=(\eps-1)E_{0,n}^-$ and let $\FF_0$ denote the far field of the local scattered electric field.
For a prescribed scalar $v$ on $\Gamma$, solve the same local Maxwell
equations with $\jump{\E_\tau}=\nabla_\Gamma v$,
$\jump{\HH_\tau}=0$, and no incident field. Define
\begin{equation}\label{eq:TL}
 \TT v=(\eps-1)E_n^-,\qquad \LL v=\FF.
\end{equation}

\begin{lemma}[Classical trace and observation maps]\label{lem:classical}
Under \eqref{eq:regime}, for fixed smooth separated particles,
\begin{equation}\label{eq:classicalmaps}
 \TT:H^{r+1}(\Gamma)\longrightarrow H^r(\Gamma),\qquad
 \LL:H^{-m}(\Gamma)\longrightarrow L_t^2(\Sph)
\end{equation}
are bounded for every real $r$ and every $m\geq0$, with constants uniform
on $I$. For rough scalar data the observation map is defined by its
continuous extension from $C^\infty(\Gamma)$; an energy solution is
asserted only for data in the natural Maxwell trace space. Furthermore,
$\norm{f_0}_{H^r}\leq C_r\norm{\mathbf a}_{L^2(\Sph)}$ for each $r$.
\end{lemma}
\begin{proof}
Appendix~\ref{app:transmission} constructs the transmission maps using
interior and exterior Maxwell admittances. The Hodge decomposition
separates the two polarizations and yields an elliptic transmission
operator of index zero. Loss gives uniqueness, and hence invertibility.
The resulting normal trace map has order one, with principal symbol
$(1-\eps)|\xi|/(1+\eps)$. The exterior representation has a smooth
kernel away from $\Gamma$, which proves the observation estimates.
All inverses depend continuously on $z\in I$.
\end{proof}

Since $\eps$ does not depend on $\diff$, neither do $\TT$, $\LL$,
$f_0$, and $\FF_0$. Diffusion enters the comparison below only through
$\NN$ and the scalar $\sg$.

Write $\mathscr F_\eta^{\gnor}(z,\diff)$ and
$\mathscr F_\eta^{\eff}(z,\diff)$ for the maps taking the Herglotz density
$\mathbf a\in L_t^2(\Sph)$ to the corresponding electric far field,
so that, for instance, $\FF_\eta^{\gnor}=\mathscr F_\eta^{\gnor}(z,\diff)\mathbf a$.

\begin{theorem}[First-order far-field equivalence]\label{thm:farfield}
Assume \eqref{eq:regime}. There are $C,\eta_0>0$ such that, for
$0<\eta<\eta_0$, the effective problem \eqref{eq:regsurface} is uniquely
solvable and
\begin{equation}\label{eq:farfielderr}
 \sup_{z\in I,\ 0\leq\diff\leq\diff_*}
 \|\mathscr F_\eta^{\gnor}(z,\diff)-\mathscr F_\eta^{\eff}(z,\diff)\|_
 {\mathcal L(L_t^2(\Sph),L_t^2(\Sph))}\leq C\eta^2.
\end{equation}
Both far fields have the expansion
\begin{equation}\label{eq:farfieldexp}
 \FF_\eta=\FF_0+\frac\eta{\sg}\LL f_0
                       +O_{L^2(\Sph)}(\eta^2\norm{\mathbf a}).
\end{equation}
\end{theorem}
\begin{proof}
Write $X=H^{-1/2}(\Gamma)$. Exact superposition gives
\begin{equation}\label{eq:traceidentity}
 f_\eta=f_0+\TT\NN f_\eta,\qquad
 \FF_\eta^{\gnor}=\FF_0+\LL\NN f_\eta.
\end{equation}
Substitute the first identity once into the second:
\begin{equation}\label{eq:twoterms}
 \FF_\eta^{\gnor}=\FF_0+\LL\NN f_0
                           +\LL\NN\TT\NN f_\eta.
\end{equation}
The middle term is $(\eta/\sg)\LL f_0+O(\eta^2\norm{\mathbf a})$ by
Lemma~\ref{lem:ntdexp}. For the last term, duality gives, for
$b\in L_t^2(\Sph)$,
\[
 (\LL\NN\TT\NN f_\eta,b)
 =\langle f_\eta,\NN^*\TT^*\NN^*\LL^*b\rangle.
\]
The map $\LL^*$ has arbitrarily smooth range, $\TT^*$ has order one,
and each application of $\NN^*$ is $O(\eta)$ between sufficiently
smooth spaces by \eqref{eq:Nsmall}. For example, start in $H^8$ and use
$H^8\to H^5\to H^4\to H^1$. Thus
\[
 \norm{\NN^*\TT^*\NN^*\LL^*b}_{H^{1/2}(\Gamma)}
 \leq C\eta^2\norm b_{L^2(\Sph)}.
\]
Corollary~\ref{cor:weaktrace} proves \eqref{eq:farfieldexp} for the
GNOR field without assuming smoothness of $f_\eta$ uniformly in
$\eta$.

For the effective model, using the equivalent Sobolev norms defined by
spectral calculus for $-\Delta_\Gamma$ yields
\[
 \norm{R_\eta}_{H^r\to H^r}\leq1,\qquad
 \norm{R_\eta}_{H^r\to H^{r+1}}\leq C\eta^{-1/2}.
\]
Since $\sg^{-1}$ is uniformly bounded, it follows that
$\norm{\TT\mathcal N_\eta^{\eff}}_{X\to X}\leq C\sqrt\eta$.
For small $\eta$ the equation
$(\Id-\TT\mathcal N_\eta^{\eff})f_\eta^{\eff}=f_0$
is invertible by a Neumann series, and
$\norm{f_\eta^{\eff}}_X\leq C\norm{\mathbf a}$.
The density determines $v=\mathcal N_\eta^{\eff}f_\eta^{\eff}\in
H^{1/2}(\Gamma)$ and hence an energy solution of the transmission problem.
Conversely every such solution satisfies this trace equation. This proves
existence and uniqueness of the effective fields.
Moreover
\[
 \mathcal N_\eta^{\eff}-\frac\eta{\sg}\Id
        =\frac{\eta^2}{\sg} R_\eta\Delta_\Gamma.
\]
It has the same smooth-data first-order estimate as \eqref{eq:Napprox}.
Apply \eqref{eq:twoterms} with $\mathcal N_\eta^{\eff}$ in place of $\NN$;
the identical duality argument proves \eqref{eq:farfieldexp} for the
effective field. Subtraction proves \eqref{eq:farfielderr}.
\end{proof}

\begin{remark}\label{rem:limits}
The constant in Theorem~\ref{thm:farfield} depends on the geometry, the
frequency interval, the damping, and $\diff_*$. The theorem is not
uniform as a gap closes, as $\gamma\downarrow0$, or for unbounded or
frequency-dependent diffusion. It concerns the regularized law
\eqref{eq:regsurface}, not unrestricted solvability of the unregularized
law $v=(\eta/\sg)f$. It approximates the exterior observation operator and
does not reconstruct the microscopic field inside the screening layer.
\end{remark}

\begin{corollary}[Second-order observation expansion]\label{cor:secondorder}
Under the hypotheses of Theorem~\ref{thm:farfield}, uniformly for
$z\in I$ and $\diff\in[0,\diff_*]$,
\begin{align}
 \FF_\eta^{\gnor}
 &=\FF_0+\frac\eta{\sg}\LL f_0
   +\frac{\eta^2}{\sg^2}\LL(\mathcal H f_0+\TT f_0)
   +O_{L^2(\Sph)}(\eta^3\|\mathbf a\|),\label{eq:secondorderhd}\\
 \FF_\eta^{\eff}
 &=\FF_0+\frac\eta{\sg}\LL f_0
   +\eta^2\LL\left(\frac1{\sg}\Delta_\Gamma f_0
                    +\frac1{\sg^2}\TT f_0\right)
   +O_{L^2(\Sph)}(\eta^3\|\mathbf a\|).
                         \label{eq:secondordereff}
\end{align}
In particular,
\begin{equation}\label{eq:regularizationdefect}
 \FF_\eta^{\gnor}-\FF_\eta^{\eff}
 =\eta^2\LL\left(\frac{\mathcal H}{\sg^2}
                  -\frac{\Delta_\Gamma}{\sg}\right)f_0
   +O_{L^2(\Sph)}(\eta^3\|\mathbf a\|).
\end{equation}
Here $\mathcal H$ denotes multiplication by mean curvature, and the
remainders are bounded operators of the incident Herglotz density.
\end{corollary}
\begin{proof}
Substitute the trace identity once more into \eqref{eq:twoterms}:
\[
 \FF_\eta^{\gnor}
 =\FF_0+\LL\NN f_0+\LL\NN\TT\NN f_0
       +\LL\NN\TT\NN\TT\NN f_\eta.
\]
For the last term, transpose the three surface maps onto $\LL^*b$.
Starting in $H^{12}$, their orders give the chain
\[
 H^{12}\xrightarrow{\NN^*}H^9
 \xrightarrow{\TT^*}H^8
 \xrightarrow{\NN^*}H^5
 \xrightarrow{\TT^*}H^4
 \xrightarrow{\NN^*}H^1.
\]
Thus this term is $O(\eta^3\|\mathbf a\|)$ by
Corollary~\ref{cor:weaktrace}. On the smooth data in the other terms,
\eqref{eq:gnorntdexp} gives
\[
 \NN=\frac\eta{\sg}\Id+\frac{\eta^2}{\sg^2}\mathcal H+O(\eta^3),
 \qquad
 \LL\NN\TT\NN f_0
 =\frac{\eta^2}{\sg^2}\LL\TT f_0+O(\eta^3\|\mathbf a\|).
\]
Every remainder here is taken between the sufficiently smooth Sobolev
spaces required by its neighboring operators, not in the natural
trace operator norm. This proves \eqref{eq:secondorderhd}.
For the regularized map, the exact identity
\[
 \mathcal N_\eta^{\eff}
 =\frac\eta{\sg}\Id+\frac{\eta^2}{\sg}\Delta_\Gamma
          +\frac{\eta^3}{\sg}R_\eta\Delta_\Gamma^2
\]
gives \eqref{eq:secondordereff} by the same argument. Subtraction proves
\eqref{eq:regularizationdefect}.
\end{proof}

The two terms at order $\eta^2$ in \eqref{eq:secondorderhd} describe
curvature and the change in the self-consistent Maxwell field,
respectively. Formula \eqref{eq:regularizationdefect} identifies the
second-order error introduced by the chosen regularization. A
second-order observation formula is therefore available from the local
transmission maps. At fixed loss, the first-order observation formula
already has the $O(\eta^2)$ accuracy of the first-order surface model.
The distinction is relevant near a simple pole. If a pole has the expansion
$z_0+\eta z_1+O(\eta^2)$, the first-order expansion of its pole factor
about $z_0$ involves the ratio $\eta |z_1|/|z-z_0|$. At $z=\operatorname{Re}z_0$ the denominator equals
$|\operatorname{Im}z_0|$, the distance of the pole from the real axis.
In the retarded problem radiation contributes to this distance. Solving the coupled surface problem allows the pole to move
with $\eta$. Proposition~\ref{prop:sphere} quantifies the shift for
fixed quasistatic spherical multipoles, and Figure~\ref{fig:resonantcmp}
illustrates the distinction for a retarded dipole.

We now match the second-order observation expansion by a regularized
surface problem.

\begin{theorem}[Second-order far-field equivalence]\label{prop:secondordersurface}
Assume \eqref{eq:regime}, and let $R_\eta$ be defined by
\eqref{eq:regN}. Set
\begin{equation}\label{eq:secondordersurfacemap}
 \NNtwo
 =\frac\eta{\sg}(2R_\eta-R_\eta^2)
   +\frac{\eta^2}{\sg^2}R_\eta\mathcal H R_\eta,
\end{equation}
where $\mathcal H$ is multiplication by mean curvature.
Replace the scalar equation in \eqref{eq:regsurface} by
$v=\NNtwo f$, with $f=(\eps-1)E_n^-$, and retain the other
transmission conditions and the local bulk Maxwell equations.
For sufficiently small $\eta$, this problem has a unique solution.
If $\widehat{\mathscr F}^{(2)}_\eta(z,\diff)$ is its electric far-field
operator, then
\begin{equation}\label{eq:secondordersurfaceerror}
 \sup_{z\in I,\ 0\leq\diff\leq\diff_*}
 \big\|\mathscr F_\eta^{\gnor}(z,\diff)
       -\widehat{\mathscr F}^{(2)}_\eta(z,\diff)\big\|
       _{\mathcal L(L_t^2(\Sph),L_t^2(\Sph))}
 \leq C\eta^3.
\end{equation}
The constants depend on the fixed geometry, frequency interval,
damping, and $\diff_*$.
\end{theorem}
\begin{proof}
Put $X=H^{-1/2}(\Gamma)$. The bounds for $R_\eta$ used in
Theorem~\ref{thm:farfield}, together with boundedness of multiplication
by $\mathcal H$ on every Sobolev space, give
\[
 \|\NNtwo\|_{X\to H^{1/2}(\Gamma)}
 \leq C(\eta^{1/2}+\eta^{3/2}),
 \qquad
 \|\TT\NNtwo\|_{X\to X}\leq C\eta^{1/2}.
\]
Hence
$\widehat f^{(2)}_\eta=f_0+\TT\NNtwo\widehat f^{(2)}_\eta$
has a unique solution with
$\|\widehat f^{(2)}_\eta\|_X\leq C\|\mathbf a\|$.
The trace reconstruction in Lemma~\ref{lem:classical} proves existence
and uniqueness of the effective Maxwell fields.

The identities
\begin{align*}
 2R_\eta-R_\eta^2
   &=\Id-\eta^2R_\eta^2\Delta_\Gamma^2,\\
 R_\eta\mathcal H R_\eta-\mathcal H
   &=\eta R_\eta\Delta_\Gamma\mathcal H R_\eta
        +\eta\mathcal H R_\eta\Delta_\Gamma
\end{align*}
show that, for every integer $r\geq0$,
\begin{equation}\label{eq:secondordersurfacesmooth}
 \left\|\left(\NNtwo-\frac\eta{\sg}\Id
              -\frac{\eta^2}{\sg^2}\mathcal H\right)f\right\|_{H^r}
 \leq C_r\eta^3\|f\|_{H^{r+4}}.
\end{equation}

The same estimate holds for the adjoint, with the scalar coefficients
conjugated, since $R_\eta$ and $\mathcal H$ are self-adjoint.
Also, $\NNtwo$ and its adjoint are $O(\eta)$ from $H^{r+4}$ to $H^r$.

Write $\widehat{\FF}^{(2)}_\eta
=\widehat{\mathscr F}^{(2)}_\eta(z,\diff)\mathbf a$ for the corresponding
far field. Iterating the coupled trace equation twice gives
\begin{align*}
 \widehat{\FF}^{(2)}_\eta
 ={}&\FF_0+\LL\NNtwo f_0
     +\LL\NNtwo\TT\NNtwo f_0\\
    &+\LL\NNtwo\TT\NNtwo\TT\NNtwo\widehat f^{(2)}_\eta.
\end{align*}
For $b\in L_t^2(\Sph)$, apply the adjoints in the last term to $\LL^*b$ as in
Corollary~\ref{cor:secondorder}. Starting in $H^{15}$ gives the chain
$H^{15}\to H^{11}\to H^{10}\to H^6\to H^5\to H^1$,
with a factor $\eta$ from each of the three surface maps.
The uniform $X$ bound on $\widehat f^{(2)}_\eta$ therefore makes this
term $O_{L^2}(\eta^3\|\mathbf a\|)$.
Applying \eqref{eq:secondordersurfacesmooth} to the smooth data in
the other terms yields
\[
 \widehat{\FF}^{(2)}_\eta
 =\FF_0+\frac\eta{\sg}\LL f_0
   +\frac{\eta^2}{\sg^2}\LL(\mathcal H f_0+\TT f_0)
   +O_{L^2(\Sph)}(\eta^3\|\mathbf a\|).
\]
Comparison with \eqref{eq:secondorderhd} proves
\eqref{eq:secondordersurfaceerror}.
\end{proof}

The second-order model can be implemented with one additional scalar
surface unknown $u$. Its scalar equations are
\begin{equation}\label{eq:secondordersurfacelocal}
 (\Id-\eta\Delta_\Gamma)u=f,\qquad
 (\Id-\eta\Delta_\Gamma)v
   =\frac\eta{\sg}(2f-u)+\frac{\eta^2}{\sg^2}\mathcal H u,
 \qquad f=(\eps-1)E_n^-.
\end{equation}
Eliminating $u$ gives \eqref{eq:secondordersurfacemap}.
As for the first-order model, this is a mathematical regularization;
no microscopic high-wave-number law or passive realization at all
frequencies is asserted.

\begin{corollary}[Exterior observations and plane waves]\label{cor:exterior}
Let $K$ be a fixed bounded smooth observation domain with
$\overline K\subset\R^3\setminus\overline D$, and let $m\geq0$ be an
integer. Under \eqref{eq:regime},
\begin{equation}\label{eq:exteriorerror}
 \sup_{z\in I,\ 0\leq\diff\leq\diff_*}
 \|\E_\eta^{\gnor}-\E_\eta^{\eff}\|_{H^m(K)}
 \leq C_{K,m}\eta^2\|\mathbf a\|_{L^2(\Sph)}.
\end{equation}
If $\widehat{\E}^{(2)}_\eta$ denotes the electric field of the
second-order surface model, then
\begin{equation}\label{eq:secondorderexteriorerror}
 \sup_{z\in I,\ 0\leq\diff\leq\diff_*}
 \big\|\E_\eta^{\gnor}-\widehat{\E}^{(2)}_\eta\big\|_{H^m(K)}
 \leq C_{K,m}\eta^3\|\mathbf a\|_{L^2(\Sph)}.
\end{equation}
The constants depend on the distance of $K$ from $\Gamma$.
Theorems and expansions above also hold for an incident plane wave
$\E^{\inc}=\mathbf p e^{\ii\kappa d\cdot x}$, $d\in\Sph$,
$d\cdot\mathbf p=0$, with $\|\mathbf a\|$ replaced by $|\mathbf p|$;
the constants are uniform in $d$.
\end{corollary}
\begin{proof}
The map taking $v$ to the electric transmission field on $K$ is
smoothing, by Appendix~\ref{app:transmission}. Replace $\LL$ by this
map in \eqref{eq:twoterms} and repeat the duality argument. The same
replacement gives the second-order expansion on $K$.
Applying it in the proof of Theorem~\ref{prop:secondordersurface}
gives \eqref{eq:secondorderexteriorerror}.
A plane wave is not an $L^2$ Herglotz superposition. Nevertheless, its
derivatives of every fixed order on a bounded set are bounded by
$C|\mathbf p|$, uniformly in $d$ and $z\in I$. These bounds provide
the incident-field estimates used in Lemmas~\ref{lem:energy} and
\ref{lem:classical}, so the proofs apply directly.
\end{proof}

\section{Resonances and a spherical verification}\label{sec:resonance}

A resonance is a pole of a meromorphically continued outgoing solution
operator in complex $z$. An observable scattering pole requires a
nonzero coupling to the chosen incident and observation channels.
A dark mode, whose coupling to these channels vanishes, need not appear
in a given far-field matrix element. Pole comparisons are based on analytic
perturbation of simple characteristic values. See \cite{GohbergSigal}
and \cite[Chapter~1]{AmmariBook}. Real-frequency error estimates alone
do not imply pole estimates. For the quasistatic hydrodynamic spectral
problem and its modal expansions on Lipschitz domains,
see~\cite{LeeRuizYu}.

For a unit sphere, a separate quasistatic calculation verifies the
surface reduction mode by mode. This calculation concerns the electrostatic
model, obtained by neglecting retardation; its poles are not the full
retarded Maxwell poles. Let $i_l$ denote the regular modified spherical Bessel function,
with prime denoting differentiation of its argument, and let $Y_{lm}$
be a spherical harmonic of degree $l\geq1$, $-l\leq m\leq l$.
Put $r_l=(1+\eta l(l+1))^{-1}$. Separation of variables gives
\begin{equation}\label{eq:sphN}
 \NN Y_{lm}=\nu_l(z,\eta)Y_{lm},\qquad
 \nu_l=\frac{i_l(\Qg)}{\Qg i_l'(\Qg)}.
\end{equation}
Since the mean curvature is $1$ on the unit sphere, the first-order and
second-order surface eigenvalues are
\begin{equation}\label{eq:gnorsphereN}
 \nu_l^{\eff}=\frac\eta{\sg}r_l,\qquad
 \widehat\nu_l^{(2)}
   =\frac\eta{\sg}(2r_l-r_l^2)
       +\frac{\eta^2}{\sg^2}r_l^2.
\end{equation}
For fixed $l$, Lemma~\ref{lem:ntdexp} gives
$\nu_l=\eta/\sg+\eta^2/\sg^2+O_l(\eta^3)$, uniformly in a complex
neighborhood on which $\varg\neq0$ and $\sg^2$ lies in a closed right
half-plane sector bounded away from zero. In particular, the scalar
Neumann problem has no pole there. The reciprocal of this
Neumann-to-Dirichlet eigenvalue appears in the multipole-dependent
coating coefficient of~\cite{LuoLAM}.

The local resonance equation $l\eps(z)+l+1=0$ has the roots
\[
 z_l^{0,\pm}=-\frac{\ii\gamma}{2}
       \pm\sqrt{\frac{l}{2l+1}-\frac{\gamma^2}{4}}.
\]
They are simple unless $\gamma^2=4l/(2l+1)$, when the two roots
coalesce. The simple-root hypothesis below excludes this case.
These roots do not depend on $\diff$.

\begin{proposition}[Quasistatic spherical poles and residues]\label{prop:sphere}
Fix $l\geq1$, $\diff\geq0$, and a simple root $z_l^0$ of $l\eps+l+1=0$.
On a sufficiently small complex neighborhood of $z_l^0$ on which $z$,
$\Omega$, $\eps$, and $\varg$ are nonzero and $\sg^2$ lies in a closed
right half-plane sector bounded away from zero, the GNOR resonance
equation in degree $l$ is
\begin{equation}\label{eq:sphres}
 l\eps+(l+1)+l(l+1)(\eps-1)\nu_l=0.
\end{equation}
The continuations of $z_l^0$ for the GNOR problem and for the first-order
and second-order surface models satisfy
\begin{equation}\label{eq:spherror}
 |z_l^{\gnor}-z_l^{\eff}|\leq C_l\eta^2,
 \qquad
 |z_l^{\gnor}-\widehat z_l^{(2)}|\leq C_l\eta^3.
\end{equation}
The roots are the unique continuations of the chosen local root in a
fixed small disc. Residues of the scalar scattering coefficient in that
multipole channel differ by the same orders. Moreover,
\begin{equation}\label{eq:gnorshift}
 z_l^{\gnor}
 =z_l^0+
 \eta\,\frac{l(l+1)}
 {(2l+1)(2z_l^0+\ii\gamma)\sg(z_l^0)}
 +O(\eta^2).
\end{equation}
In the separate case $\gamma=0$, with $z_l^0=\sqrt{l/(2l+1)}$, this
becomes
\begin{equation}\label{eq:gnorzeroloss}
 z_l^{\gnor}
 =z_l^0+\frac\eta2\sqrt{l(l+1)}
                  \sqrt{1-\ii\diff z_l^0}+O(\eta^2).
\end{equation}
At $\diff=0$ this is the hydrodynamic blueshift
$z_l^0+\frac\eta2\sqrt{l(l+1)}+O(\eta^2)$. For $\diff>0$ the first
correction has positive real part and negative imaginary part.
\end{proposition}
\begin{proof}
Let $r=|x|$ be the radial coordinate, and use $Y_{lm}$ as shorthand
for $Y_{lm}(x/|x|)$. For a source-free electrostatic mode take
$\phi^-=Ar^lY_{lm}$ and $\phi^+=Br^{-l-1}Y_{lm}$, with
$\ET=-\nabla\phi^-$ and $\E^+=-\nabla\phi^+$.
The normal flux condition gives $(l+1)B=-\eps lA$.
The tangential condition gives
$B-A=l(\eps-1)\nu_l A$.
Eliminating $A,B$ yields \eqref{eq:sphres}.
Replacing $\nu_l$ by $\nu_l^{\eff}$ or $\widehat\nu_l^{(2)}$ gives the
effective equations. For either comparison, let $D_{\gnor}$ and
$D_{\eff}$ denote the two left-hand sides, and set $p=2$ or $3$,
respectively. Choose a positively oriented circle $C$ about $z_l^0$
whose closed disc lies in the stated neighborhood and contains no other
local root. Both denominators converge uniformly there to $l\eps+l+1$.
By Rouch\'e's theorem each has a unique simple zero in the disc, denoted
by $z_{\gnor}$ and $z_{\eff}$. The expansion of $\nu_l$ gives
$D_{\gnor}-D_{\eff}=O(\eta^p)$ on a slightly larger disc. By Cauchy's
estimate the same bound holds for the derivative difference on $C$.
Both denominators stay uniformly away from zero on $C$ and their
derivatives remain bounded there. The argument principle yields
\[
 z_{\gnor}-z_{\eff}
 =\frac{1}{2\pi\ii}\int_C z\left(
 \frac{D_{\gnor}'}{D_{\gnor}}-\frac{D_{\eff}'}{D_{\eff}}\right)\dd z
 =O(\eta^p).
\]
This proves \eqref{eq:spherror}. To compare residues,
take $\phi^+=c_lr^lY_{lm}+d_lr^{-l-1}Y_{lm}$. The same two interface
conditions give the scalar scattering coefficient
\begin{equation}\label{eq:scalarmultipole}
 S_l(z,\eta):=\frac{d_l}{c_l}
 =-\frac{l(\eps-1)(1-l\nu_l)}
 {l\eps+l+1+l(l+1)(\eps-1)\nu_l}.
\end{equation}
Replacing $\nu_l$ by a surface eigenvalue changes both numerator and
denominator by $O(\eta^p)$ on $C$. The two scattering coefficients thus
differ by $O(\eta^p)$ there, and integration gives the residue estimate.
At $\eta=0$ the numerator equals $2l+1$ at $z_l^0$.
At the local root $\Omega(z_l^0)=l/(2l+1)$, and
$\eps'(z)=\Omega'(z)/\Omega(z)^2$. Differentiating the local
denominator and using the leading term $\eta/\sg$ gives
\eqref{eq:gnorshift}. When $\gamma=0$,
$\sg(z_l^0)=\sqrt{(l+1)/(2l+1)}/\sqrt{1-\ii\diff z_l^0}$,
which proves \eqref{eq:gnorzeroloss}.
The simplicity here is that of the scalar radial equation; the usual
$2l+1$ angular degeneracy is unchanged.
\end{proof}

The negative imaginary correction gives diffusion-induced damping even
when the Drude collision rate is zero. This shift and broadening
mechanism is already part of GNOR \cite{GNOR,RazaReview}. The estimates
above compare the exact and regularized surface poles. They concern
fixed quasistatic multipoles. The $\gamma=0$ calculation is an analytic
scalar calculation, not an extension of the uniform Maxwell theorem
to vanishing collision loss.

For the dipole, let $\mathbf p_{\mathrm{dip}}$ be the induced
electric dipole moment in units of $\eps_0a^3E_*$, and define the
polarizability $\alpha$ by $\mathbf p_{\mathrm{dip}}=\alpha\E^{\inc}$.
In the same quasistatic model,
\begin{equation}\label{eq:polarizability}
 \frac\alpha{4\pi}=\frac{\eps-1-t}{\eps+2+2t},\qquad
 t=(\eps-1)\nu_1.
\end{equation}
More generally the scalar local permittivity reproducing the exterior
degree-$l$ quasistatic scattering coefficient is
\begin{equation}\label{eq:fit}
 \eps_{\mathrm{fit},l}
       =\frac{\eps}{1+l(\eps-1)\nu_l},
\end{equation}
away from zeros of this denominator. It depends on $l$.
Thus even for a sphere there is generally no single scalar bulk parameter
which reproduces every multipole.

\section{Concluding remarks}\label{sec:conclusion}

We have established quantitative surface approximations for Maxwell
scattering by smooth, separated metallic particles described by the
hard-wall GNOR model, which contains the hydrodynamic Drude model as the
case $\diff=0$. The length $\ell=\beta/\omega_p$ sets the reference
screening scale. For a planar mode at zero tangential wave number, the
physical decay length is $a/\operatorname{Re}\Qg=\ell/\operatorname{Re}\sg$.
The ratio $\eta=\ell/a$ is the expansion parameter. A surface
approximation replaces the layer by a transmission condition on the
particle surface, so that only local Maxwell equations are solved in
the bulk. The first-order and second-order conditions reproduce the
scattered far field of the nonlocal problem with errors $O(\eta^2)$ and
$O(\eta^3)$ (Theorems~\ref{thm:farfield}
and~\ref{prop:secondordersurface}). The same bounds hold for
the scattered field at any fixed distance from the particles. They are
uniform on compact subplasma frequency intervals, for fixed positive
damping and bounded diffusion. They concern the field outside the
particles, not the field inside the screening layer.

The argument rests on three observations. The electric field in the
metal splits exactly into a divergence-free part, which obeys the local
Drude equations, and the gradient of a scalar potential, which decays
exponentially across the screening layer
(Proposition~\ref{prop:exact}). Eliminating the potential leaves a
local problem whose surface operator is $O(\eta)$ on smooth data in the
sense of \eqref{eq:Nsmall}. Its second-order term contains the mean
curvature. The term $\TT f_0$ in \eqref{eq:secondorderhd} arises from
solving the coupled Maxwell problem. The proof requires only the uniform
$H^{-1/2}(\Gamma)$ bound for the exact boundary density, since the
far-field map is smoothing. Diffusion enters through the current
dissipation and the screening parameter $\sg$.

The surface condition can also be realized by a thin layer of local
material, which suits standard electromagnetic solvers
(Appendix~\ref{app:coating}). The anisotropic vacuum-side coating may
have physical thickness comparable to $\ell$. For the isotropic metal-side
coating, physical thickness $O(a\eta^2)$ suffices for the same
$O(\eta^2)$ accuracy. These representations are justified for a
planar interface at bounded tangential wave numbers and for a sphere on
finitely many multipoles. At $\diff=0$ the vacuum-side layer is not
passive, so it is a computational device rather than a physical
material. A retarded dipole computation on a sphere confirms all
predicted rates for several values of $\diff$.

Table~\ref{tab:scope} summarizes the scope of the estimates. Resonance
poles and residues are compared only for fixed quasistatic multipoles of
a sphere (Proposition~\ref{prop:sphere}). Extending this comparison to
an isolated retarded resonance of a nonspherical particle requires
uniform estimates on a complex-frequency contour enclosing the pole. Further questions lie beyond the
fixed-geometry setting. They include estimates uniform as two particles
approach contact, where the gap width becomes a second small parameter,
and the derivation of a quantitative surface law from a self-consistent
microscopic model with spill-out or tunnelling, for which time-dependent density-functional theory calculations offer one starting point \cite{Bonafe}. They also include
causal and passive effective models, particles with corners, and
unbounded diffusion (Remark~\ref{rem:unboundeddiff}).

\begin{table}[htbp]
\centering\small
\begin{tabular}{@{}>{\raggedright\arraybackslash}p{.24\textwidth}
 >{\raggedright\arraybackslash}p{.31\textwidth}
 >{\raggedright\arraybackslash}p{.38\textwidth}@{}}
\toprule
Setting & Quantitative conclusion & Scope\\
\midrule
Fixed smooth particles & $O(\eta^2)$ first-order and $O(\eta^3)$
second-order surface-model errors & Fixed loss and subplasma interval;
all $L^2$ Herglotz densities; uniform for $0\leq\diff\leq\diff_*$\\[4pt]
Thin planar coatings & $O(h\eta+\eta^3)$ anisotropic error;
$O(h+\eta^3)$ isotropic error & Bounded tangential wave numbers;
scattering estimates away from grazing incidence\\[4pt]
Coated sphere & $O(\eta^2+h\eta)$ anisotropic error;
$O(\eta^2+h)$ isotropic error & Retarded Maxwell problem on each fixed
finite angular subspace\\[4pt]
Quasistatic sphere & $O(\eta^2)$ first-order and $O(\eta^3)$
second-order pole and residue comparison & Fixed
multipole, fixed $\diff\geq0$, and simple radial root; stated complex
neighborhood\\
\bottomrule
\end{tabular}
\caption{Regimes and conclusions. Constants depend on the fixed
frequency window and geometry or cutoff in the corresponding statement. The isotropic coating
estimates assume $0<h\leq c_*\eta$ with $c_*$ sufficiently small.}
\label{tab:scope}
\end{table}

\section*{Statements and declarations}

\paragraph{Funding.}

H.L. was supported by the National Research Foundation of Korea (NRF) grant No. NRF-2018R1D1A1B07042678 and grant No. NRF-2022R1A4A5033271. S.Y. was supported by the National Research Foundation of Korea (NRF) grant No. RS-2025-24523302 and from Korea University grant No. K1923451.

\paragraph{Competing interests.}
The authors have no competing interests to declare.

\paragraph{Data availability.}
The Python scripts supplied with the manuscript reproduce the numerical
figures and tables in Appendix~\ref{app:coating}. They use NumPy,
SciPy, and mpmath. All data are generated by the spherical matching
formulas.

\subsection*{Declaration of AI use}
The authors interacted with and guided ChatGPT~6 and Claude Opus~5.5 whilst developing the proofs and writing the paper. The authors take full responsibility for the accuracy and integrity of all content in this paper.

\appendix
\section{Thin-layer realization: quantitative estimates}\label{app:coating}

We justify the coating construction as a mathematical representation
of the first-order surface law in two settings: a planar interface
on bounded tangential frequency bands, and a sphere on any fixed finite
set of Maxwell multipoles. Retardation is retained. The estimates are
uniform in the thin-layer contrast, which may diverge or degenerate.
They do not assert a uniform approximation over all tangential frequencies
or all multipoles. Thin-layer transmission conditions have an established
asymptotic theory; see, for example,~\cite{PeronThin}. The estimates in the following specifically address the coefficients used in the local analogue
of~\cite{LuoLAM} and in the surface-moment construction.
Related exact planar scattering formulas are given in \cite{Pitelet}.
Passivity of the auxiliary coating as a standalone material is not
assumed; see the final subsection of this appendix.

Throughout this appendix, $z\in I$ and $0\leq\diff\leq\diff_*$, with
$I$, $\gamma$, $\sigma$, $s$, $\sg$, $\Qg$, and $\eps$ as in
Section~\ref{sec:notation}. Let $h>0$ denote the constant
dimensionless coating thickness and set
\begin{equation}\label{eq:coatb}
 b_0(z,\eta)=(1-\eps^{-1})\Qg^{-1}=\frac{\eta}{s^2\sg}
       =\frac{\eta\sqrt{\varg}}{s^3}.
\end{equation}
Both $s$ and $\sg$ occur here because $1-\eps^{-1}=s^{-2}$ is a
transverse identity. At $\diff=0$, $b_0=\eta/s^3$.
Thus, $c\eta\leq |b_0|\leq C\eta$ uniformly on $I\times[0,\diff_*]$.
For real $z\in I$, $\arg s\in(-\pi/4,0)$ and
$\arg\sqrt{\varg}\in(-\pi/4,0]$, so $\arg b_0\in(-\pi/4,3\pi/4)$.
Compactness gives $\operatorname{Re}(e^{-\ii\pi/4}b_0)\geq c|b_0|$,
and therefore
\begin{equation}\label{eq:coatnoncancel}
 |h+b_0|\geq c(h+|b_0|),\qquad h>0.
\end{equation}
All permeabilities are $1$. Constants may depend on $I$, $\gamma$,
$\sigma$, $\diff_*$, and a stated frequency or multipole cutoff, but
not on $h$, $\eta$, or $\diff$.

\subsection{Matching coefficients and the reference surface}

Consider metal in $t<0$ and vacuum in $t>0$. At zero tangential wave
number, the electrostatic GNOR surface law is
$\phi^+-V=b_0\partial_t\phi^+$, where $V$ is the effective metal
potential extrapolated to $t=0$. Replace the metal strip $-h<t<0$
by an isotropic coating of permittivity $\eps_c$. Let $\phi_c$ and
$\phi_m$ denote the coating and local metal potentials. Continuity of normal
displacement gives the constant
$q=\partial_t\phi^+=\eps_c\partial_t\phi_c
=\eps\partial_t\phi_m$. Hence, exactly,
\begin{equation}\label{eq:metallayer}
 \phi^+(0)-V(0)=h\left(\frac1{\eps_c}-\frac1\eps\right)q,
 \qquad \eps_c=\frac{h}{b_0+h/\eps}.
\end{equation}
The subtraction of $h/\eps$ accounts for the removed metal. In particular,
at fixed $b_0\ne0$,
\begin{equation}\label{eq:localanalogue}
 \lim_{h\downarrow0}\frac{\eps_c}{h}
 =\frac1{b_0}=\frac{\eps\Qg}{\eps-1}.
\end{equation}
At $\diff=0$ this is the planar local-analogue coefficient
of~\cite{LuoLAM} in the present normalization. Finite-frequency estimates require more than this
normal-flux identity.

\begin{figure}[htbp]
\centering
\begin{tikzpicture}[x=0.94cm,y=0.94cm,font=\small]
 \begin{scope}
  \node[font=\small\bfseries] at (3.15,4.0) {(a) Anisotropic vacuum-side coating};
  \fill[metalwash] (0.2,1.1) rectangle (3.0,3.0);
  \fill[screentint] (3.0,1.1) rectangle (3.8,3.0);
  \draw[figboundary] (3.0,1.1)--(3.0,3.0);
  \draw[screenorange,thick] (3.8,1.1)--(3.8,3.0);
  \node at (1.55,2.15) {metal, $\eps$};
  \node[align=center] at (3.4,2.15) {$\eps_\tau^c$\\$\eps_n^c$};
  \node at (5.15,2.15) {vacuum};
  \draw[figarrow] (0.2,0.9)--(6.2,0.9) node[right] {$t$};
  \draw (3.0,0.82)--(3.0,0.98) node[below=5pt] {$0$};
  \draw (3.8,0.82)--(3.8,0.98) node[below=5pt] {$h$};
  \draw[figarrow,fieldteal] (5.4,3.35)--(3.0,3.35);
  \node[font=\scriptsize,fieldteal,above] at (4.2,3.35) {extrapolate through vacuum};
  \node[font=\scriptsize,align=center] at (3.15,0.05)
     {$\eps_\tau^c=1$, $\eps_n^c=h/(h+b_0)$\\$h=O(\eta)$ for $O(\eta^2)$ accuracy};
 \end{scope}
 \begin{scope}[xshift=7.25cm]
  \node[font=\small\bfseries] at (3.15,4.0) {(b) Isotropic metal-side coating};
  \fill[metalwash] (0.2,1.1) rectangle (2.2,3.0);
  \fill[screentint] (2.2,1.1) rectangle (3.0,3.0);
  \draw[screenorange,thick] (2.2,1.1)--(2.2,3.0);
  \draw[figboundary] (3.0,1.1)--(3.0,3.0);
  \node at (1.15,2.15) {metal, $\eps$};
  \node at (2.6,2.15) {$\eps_c$};
  \node at (4.6,2.15) {vacuum};
  \draw[figarrow] (0.2,0.9)--(6.2,0.9) node[right] {$t$};
  \draw (2.2,0.82)--(2.2,0.98) node[below=5pt] {$-h$};
  \draw (3.0,0.82)--(3.0,0.98) node[below=5pt] {$0$};
  \draw[figarrow,fieldteal] (0.6,3.35)--(3.0,3.35);
  \node[font=\scriptsize,fieldteal,above] at (1.8,3.35) {extrapolate through metal};
  \node[font=\scriptsize,align=center] at (3.15,0.05)
     {$\eps_c=h/(b_0+h/\eps)$\\$h=O(\eta^2)$ for $O(\eta^2)$ accuracy};
 \end{scope}
\end{tikzpicture}
\caption{Two thin-coating representations of the hard-wall surface law.
The original reference interface is $t=0$ in both panels; the normal points
toward increasing $t$. Orange regions show the auxiliary coatings.
In (a) the exterior field is extrapolated through vacuum to $t=0$;
in (b) the remaining local metal field is extrapolated through metal to
the same reference plane. The displayed thickness choices yield
$O(\eta^2)$ errors in planar transfer matrices on bounded tangential
wave-number bands, in planar scattering on the propagating bands of
Corollary~\ref{cor:coatplanarscattering}, and in scattering on fixed
finite spherical angular subspaces.
The coatings need not be passive as standalone materials.}
\label{fig:coatings}
\end{figure}
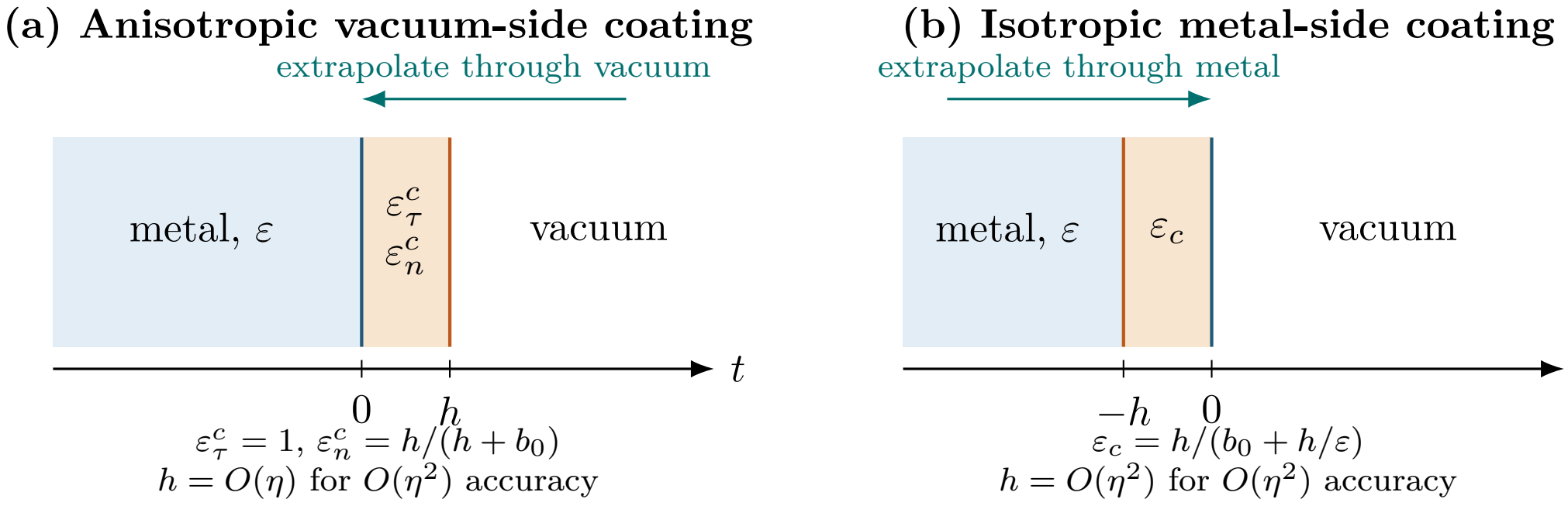

Alternatively, put a uniaxial coating in $0<t<h$, with tensor
$\eps_\tau^c(\Id-\nn\nn^\top)+\eps_n^c\nn\nn^\top$. Define its excess
tangential polarization and normal resistance by
\begin{equation}\label{eq:coatingmoments}
 p_c=h(\eps_\tau^c-1),\qquad
 b_c=h\bigl((\eps_n^c)^{-1}-1\bigr).
\end{equation}
The first-order surface moments are matched by
\begin{equation}\label{eq:coatingeps}
 p_c=d_\parallel(\eps-1),\qquad
 b_c=-d_\perp(1-\eps^{-1}),\qquad
 \eps_\tau^c=1+p_c/h,\quad \eps_n^c=(1+b_c/h)^{-1}.
\end{equation}
Here, the exterior field is extrapolated through vacuum from $t=h$ back
to $t=0$. For the hard-wall first-order law,
\begin{equation}\label{eq:anisocoat}
 p_c=0,\qquad b_c=b_0,\qquad
 \eps_\tau^c=1,\qquad \eps_n^c=\frac{h}{h+b_0}.
\end{equation}
By \eqref{eq:coatnoncancel}, $\eps_n^c$ is finite and nonzero.
The extrapolations above specify the reference interface and are essential
to the signs and finite-thickness corrections.
Figure~\ref{fig:coatings} displays the two placements and the fields
being extrapolated.

\subsection{Exact planar Maxwell transfer and its error}

After a tangential Fourier transform, rotate coordinates so that the
wave vector is $(k,0)$, with $k\geq0$ in Cartesian coordinates $(x,y,t)$, and $\nn=e_3$,
the third coordinate unit vector.
For transverse electric (TE) and transverse magnetic (TM) polarization,
use the tangential states
\[
 Y_{\mathrm{TE}}=(E_y,-H_x)^\top,
 \qquad Y_{\mathrm{TM}}=(E_x,H_y)^\top.
\]
In a homogeneous uniaxial medium, Maxwell's equations give $Y'=A^cY$,
where the prime here denotes differentiation in the normal coordinate $t$, and
\begin{equation}\label{eq:coatmatrices}
 A^c_{\mathrm{TE}}=\begin{pmatrix}0&\ii\kappa\\
 \ii(\kappa\eps_\tau^c-k^2/\kappa)&0\end{pmatrix},\qquad
 A^c_{\mathrm{TM}}=\begin{pmatrix}
 0&\ii(\kappa-k^2/(\kappa\eps_n^c))\\
 \ii\kappa\eps_\tau^c&0\end{pmatrix}.
\end{equation}
For example, in the TM channel,
$E_n=-kH_y/(\kappa\eps_n^c)$; substitution into the tangential curl
equation gives the upper-right entry. Let $A^0$ denote the corresponding
vacuum matrix. With the extrapolated exterior state
$Y^+(0)=e^{-hA^0}Y(h)$, the exact coating transfer is
\begin{equation}\label{eq:coattransfer}
 Y^+(0)=T_hY^-(0),\qquad
 T_h=e^{-hA^0}e^{hA^c}=e^{-hA^0}e^{hA^0+K},
\end{equation}
where
\begin{equation}\label{eq:coatK}
 K_{\mathrm{TE}}=\begin{pmatrix}0&0\\\ii\kappa p_c&0\end{pmatrix},
 \qquad
 K_{\mathrm{TM}}=\begin{pmatrix}0&-\ii k^2b_c/\kappa\\
                          \ii\kappa p_c&0\end{pmatrix}.
\end{equation}

\begin{lemma}[A matrix estimate]\label{lem:coatmatrix}
For matrices $A,K$ of fixed dimension, with $\|A\|+\|K\|\leq M$ and
$0<h\leq1$,
\begin{equation}\label{eq:coatmatrixbound}
 \|e^{-hA}e^{hA+K}-e^K\|
 +\|e^{hA+K}e^{-hA}-e^K\|\leq C_Mh\|K\|.
\end{equation}
Consequently, each product differs from $\Id+K$ by at most
$C_M(h\|K\|+\|K\|^2)$. If $K^2=0$, the second term is absent.
\end{lemma}
\begin{proof}
Put $W(t)=e^{-thA}e^{t(hA+K)}$ for $0\leq t\leq1$. Then
$W'=K_h(t)W$, $W(0)=\Id$, where
$K_h(t)=e^{-thA}Ke^{thA}$ and
$\|K_h(t)-K\|\leq C_Mh\|K\|$.
Variation of constants, together with
$\|W(t)\|+\|e^{tK}\|\leq C_M$, proves the estimate for the first
product. The second is its conjugate by $e^{hA}$; the extra difference
$e^{hA}e^Ke^{-hA}-e^K$ is also bounded by $C_Mh\|K\|$, as follows
by differentiating the conjugation and using
$[A,e^K]=\int_0^1 e^{(1-t)K}[A,K]e^{tK}\dd t$.
Here, $[A,B]=AB-BA$ is the matrix commutator.
Finally, the exponential series gives
$\|e^K-\Id-K\|\leq C_M\|K\|^2$, with equality to zero when $K^2=0$.
\end{proof}

\begin{theorem}[Planar moment approximation]\label{thm:planarcoat}
Fix $K_*>0$ and suppose that $0\leq k\leq K_*$,
$|p_c|+|b_c|\leq C\eta$, and $h,\eta$ are sufficiently small.
Assume that the coating coefficients in \eqref{eq:coatingeps} are finite and
$\eps_n^c\ne0$. Then, in either polarization,
\begin{equation}\label{eq:coatgeneralerror}
 \sup_{z\in I,\,0\leq k\leq K_*}
 \|T_h-(\Id+K)\|\leq C_{K_*}(h\eta+\eta^2).
\end{equation}
For the hard-wall coating \eqref{eq:anisocoat}, the sharper estimate is
\begin{equation}\label{eq:coatpureerror}
 \|T_h-(\Id+K)\|\leq C_{K_*}h\eta.
\end{equation}
The TE transfer is exactly $\Id$ in this case. These are also operator
estimates on $L^2$ tangential states whose Fourier transforms are
supported in $|\xi|\leq K_*$.
\end{theorem}
\begin{proof}
The vacuum matrices are uniformly bounded on the stated compact set,
since $\kappa=\sigma z$ is bounded away from zero. Equation
\eqref{eq:coatK} gives $\|K\|\leq C_{K_*}\eta$.
Apply Lemma~\ref{lem:coatmatrix}; for \eqref{eq:anisocoat}, $K^2=0$.
The operator statement follows from Plancherel's theorem; rotations of
the tangential coordinates preserve the norm.
\end{proof}

The law $Y^+=(\Id+K)Y^-$ is the first-order form of
\eqref{eq:feibelman}. For example, its TM electric jump is
$[E_x]=-\ii k^2b_cH_y^-/\kappa$, and its tangential magnetic jump is
$[H_y]=\ii\kappa p_cE_x^-$. Replacing zeroth-order traces by either
side in the first-order surface constitutive law changes terms only at
second order. For the hard-wall law, $H_y$ is continuous, and this
identification is exact.

There is a distinction between a thin layer and a weak surface response.
At fixed $p_c,b_c$, \eqref{eq:coattransfer} tends to $e^K$ as $h\to0$.
In the TM channel,
$K^2=k^2p_cb_c\Id$. Thus, matching two finite moments does not in
general give the exact zero-thickness transfer $\Id+K$.
The estimate \eqref{eq:coatgeneralerror} uses their joint smallness.

For comparison with the exact GNOR half-space, set
\begin{equation}\label{eq:coatbhdk}
 b_{\gnor}(k)=(1-\eps^{-1})(\Qg^2+k^2)^{-1/2}.
\end{equation}
The branch has a positive real part. Exact longitudinal elimination gives
$T^{\gnor}_{\mathrm{TE}}=\Id$ and
$T^{\gnor}_{\mathrm{TM}}=\Id+K_{\mathrm{TM}}(0,b_{\gnor})$.
Uniformly on bounded $k$-intervals,
\begin{equation}\label{eq:coatbhderror}
 |b_{\gnor}(k)-b_0|
 =|b_0|\left|(1+\eta^2k^2/\sg^2)^{-1/2}-1\right|
 \leq C_{K_*}\eta^3k^2.
\end{equation}
It follows that the anisotropic vacuum-side coating satisfies
\begin{equation}\label{eq:coatplanehd}
 \|T_h-T^{\gnor}\|\leq C_{K_*}(h\eta+\eta^3).
\end{equation}
The right-hand side is $O(\eta^2)$ when $h=O(\eta)$.

\begin{corollary}[Planar reflection and transmission]\label{cor:coatplanarscattering}
Fix $\theta\in(0,1)$ and restrict to
$0\leq k\leq(1-\theta)\kappa$. For sufficiently small $h,\eta$, the
anisotropic coating and the GNOR interface have unique TE and TM
scattering solutions. Their reflected electric amplitudes and transmitted transverse bulk
electric amplitudes, referred to $t=0$, differ by at most
$C_\theta(h\eta+\eta^3)$ times the incident amplitude. The same bound
holds in $L^2$ for incident angular spectra supported in this band.
The transmitted amplitude does not include the resolved longitudinal
field in the metal.
\end{corollary}
\begin{proof}
Let $q_0=(\kappa^2-k^2)^{1/2}>0$ and
$q_m=(\kappa^2\eps-k^2)^{1/2}$ with $\operatorname{Im}q_m>0$.
The exterior incoming and outgoing states are
$w_{\mathrm{in}}=(1,-y_0)^\top$ and $w_{\mathrm{out}}=(1,y_0)^\top$;
the decaying metal state is $w_m=(1,-y_m)^\top$. Here
$(y_0,y_m)=(q_0/\kappa,q_m/\kappa)$ for TE and
$(\kappa/q_0,\kappa\eps/q_m)$ for TM.
The amplitudes are normalized by the tangential electric component;
away from grazing incidence, this is equivalent to the usual electric
amplitude normalization. Denote the reflected and transmitted amplitude coefficients by
$c_{\mathrm{refl}}$ and $c_{\mathrm{trans}}$, respectively. They solve
\[
 \bigl(w_{\mathrm{out}},-Tw_m\bigr)\binom{c_{\mathrm{refl}}}{c_{\mathrm{trans}}}=-w_{\mathrm{in}}.
\]
At $T=\Id$ this matrix is invertible: the local passive half-space has
no nonzero outgoing solution without an incident wave, as follows from
the energy identity. Its inverse is uniformly bounded on the compact
set under consideration. Both transfers are $O(\eta)$ perturbations of
$\Id$. A matrix inverse identity and \eqref{eq:coatplanehd} prove the
claim. Plancherel gives the angular-spectrum bound.
\end{proof}

\subsection{The isotropic metal-side layer needs a thinner scaling}

The finite-thickness correction in \eqref{eq:metallayer} matches the
normal moment exactly but introduces a tangential excess. Let $A^m$ be
\eqref{eq:coatmatrices} with both permittivities equal to $\eps$.
For the strip $-h<t<0$, extrapolate the remaining metal to $t=0$.
The exact transfer is now given by
\begin{equation}\label{eq:metaltransfer}
 T_h^m=e^{hA^c}e^{-hA^m},\qquad
 h(A^c-A^m)=K(p_m,b_0),\qquad p_m=h(\eps_c-\eps),
\end{equation}
with the appropriate matrix from \eqref{eq:coatK}.

\begin{proposition}[Isotropic planar realization]\label{prop:metalcoat}
There is $c_*>0$ such that, for $0<h\leq c_*\eta$ and sufficiently
small $\eta$, the coefficient \eqref{eq:metallayer} is finite and
\begin{equation}\label{eq:metalcoaterror}
 \sup_{z\in I,\,0\leq k\leq K_*}
 \|T_h^m-T^{\gnor}\|\leq C_{K_*}(h+\eta^3).
\end{equation}
The same order holds for reflection and transmission on the propagating
bands of Corollary~\ref{cor:coatplanarscattering}. In particular,
$h=O(\eta^2)$ yields an $O(\eta^2)$ approximation.
\end{proposition}
\begin{proof}
Choose $c_*$ so that $|h/\eps|\leq |b_0|/2$ on $I$. Then
$|\eps_c|\leq C h/\eta$ and $|p_m|\leq Ch$.
Write $K(p_m,b_0)=K(0,b_0)+K(p_m,0)$.
The first summand is nilpotent, and the second is $O(h)$.
Lemma~\ref{lem:coatmatrix} and the Lipschitz estimate for the matrix
exponential on bounded sets give
$\|T_h^m-\Id-K(0,b_0)\|\leq Ch$.
Use \eqref{eq:coatbhderror} and then the same scattering matrix argument
as above.
\end{proof}

Thus, a thickness $h$ comparable to $\eta$ is adequate for the anisotropic
vacuum-side construction, but does not in general give first-order
Maxwell accuracy for an isotropic metal-side coating. The latter has
$p_m=O(\eta)$, generically nonzero at first order, when $h\asymp\eta$.
Its TE transfer already records this discrepancy.

\subsection{A curved, retarded verification on the sphere}

We now take $D$ to be the unit ball. Let $Y_{lm}$ be orthonormal scalar
spherical harmonics, $l\geq1$, put $\lambda_l=l(l+1)$, and define
$\Psi_{lm}=\nabla_{\Sph}Y_{lm}/\sqrt{\lambda_l}$ and
$X_{lm}=\nn\times\Psi_{lm}$, with $\nn=x/|x|$.
The gradient $\nabla_{\Sph}$ is the surface gradient on the unit sphere. For a radially uniaxial medium with tangential and normal relative
permittivities $\eps_\tau(r)$ and $\eps_n(r)$, introduce scalar radial
amplitudes $u(r)$ and $v(r)$. A TM mode can be written
\begin{equation}\label{eq:coatsphericalfields}
 \HH=\frac{u(r)}r X_{lm},\qquad
 \E=\frac{v(r)}{\ii\kappa r}\Psi_{lm}
 +\frac{\sqrt{\lambda_l}u(r)}{\ii\kappa\eps_n(r)r^2}Y_{lm}\nn.
\end{equation}
Direct substitution into the curl equations gives
\begin{equation}\label{eq:coatradialode}
 \frac{\dd}{\dd r}\binom uv
 =A_l^c(r)\binom uv,
 \qquad A_l^c(r)=\begin{pmatrix}
 0&\eps_\tau(r)\\
 \lambda_l/(\eps_n(r)r^2)-\kappa^2&0
 \end{pmatrix}.
\end{equation}
The state $(u,v)$ is continuous at an ordinary material interface.
Write $U_l^a(r_2,r_1)$ for the fundamental matrix in medium $a$,
mapping the state at $r_1$ to that at $r_2$. Superscripts $0,m,c$
mean vacuum, local Drude metal, and coating, respectively.

For the vacuum-side coating \eqref{eq:anisocoat} in $1<r<1+h$,
\begin{equation}\label{eq:coatsphereT}
 T_l^c=U_l^0(1,1+h)U_l^c(1+h,1),\qquad
 A_l^c(r)=A_l^0(r)+\frac{b_0}{h}B_l(r),\qquad
 B_l(r)=\begin{pmatrix}0&0\\\lambda_l/r^2&0\end{pmatrix}.
\end{equation}
For an isotropic coating \eqref{eq:metallayer} in $1-h<r<1$,
the corresponding transfer is
$T_l^m=U_l^c(1,1-h)U_l^m(1-h,1)$.

\begin{lemma}[Spherical transfer estimates]\label{lem:spherecoattransfer}
Fix $L\geq1$. Uniformly for $1\leq l\leq L$ and $z\in I$,
\begin{align}
 \|T_l^c-\Id-b_0B_l(1)\|&\leq C_Lh\eta,
                      &&0<h\leq h_0,\label{eq:spherecoatTerror}\\
 \|T_l^m-\Id-b_0B_l(1)\|&\leq C_Lh,
                      &&0<h\leq c_*\eta.\label{eq:spheremetalTerror}
\end{align}
Here $h_0\in(0,1/2)$ and $\eta$ are sufficiently small. The bounds
refer to the TM states in \eqref{eq:coatsphericalfields}; the TE channel
is treated separately below.
\end{lemma}
\begin{proof}
For the first estimate, set $r=1+ht$ and
\[
 W(t)=U_l^0(1,r)U_l^c(r,1),\qquad 0\leq t\leq1.
\]
Then $W(0)=\Id$, $W(1)=T_l^c$, and
$W'=b_0\widetilde B_h(t)W$, where
\[
 \widetilde B_h(t)=U_l^0(1,r)B_l(r)U_l^0(r,1),
 \qquad
 \sup_t\|\widetilde B_h(t)-B_l(1)\|\leq C_Lh.
\]
The last estimate follows from the boundedness and smoothness of
$A_l^0,B_l$ in a fixed neighborhood of $r=1$.
Variation of constants compares $W(1)$ with $e^{b_0B_l(1)}$ with error
$C_Lh|b_0|$. Since $B_l(1)^2=0$, this exponential is
$\Id+b_0B_l(1)$.

For the metal-side layer, put $r_0=1-h$, $r=r_0+ht$, and
\[
 W(t)=U_l^m(1,r)U_l^c(r,r_0)U_l^m(r_0,1).
\]
Again $W(0)=\Id$ and $W(1)=T_l^m$. Its generator is the conjugate
by $U_l^m(1,r)$ of
\[
 h(A_l^c-A_l^m)=
 \begin{pmatrix}0&p_m\\\lambda_lb_0/r^2&0\end{pmatrix}.
\]
This differs from $b_0B_l(1)$ by at most
$C_L(|p_m|+h|b_0|)\leq C_Lh$.
The same variation-of-constants argument proves the second estimate.
\end{proof}

Use $\nu_l(z,\eta)=i_l(\Qg)/(\Qg i_l'(\Qg))$ as in \eqref{eq:sphN}, where
$i_l$ is the regular modified spherical Bessel function. The exact GNOR surface law on the
sphere is
\begin{equation}\label{eq:coatspherehd}
 T_l^{\gnor}=\Id+b_lB_l(1),\qquad
 b_l=(1-\eps^{-1})\nu_l(z,\eta),\qquad
 b_l=b_0+O_L(\eta^2).
\end{equation}
Indeed, \eqref{eq:coatsphericalfields} gives
$E_n^- =\sqrt{\lambda_l}u/(\ii\kappa\eps)$ at $r=1$.
The exact electric jump \eqref{eq:exactbc} therefore gives
$[v]=\lambda_l b_lu$, while $[u]=0$.
The expansion follows from Lemma~\ref{lem:ntdexp} on each spherical
harmonic. For the regularized law \eqref{eq:regsurface}, the same formula
holds with
\begin{equation}\label{eq:coatspherereg}
 b_l^{\eff}=\frac{b_0}{1+\eta\lambda_l}=b_0+O_L(\eta^2).
\end{equation}

\begin{theorem}[Coated-sphere far-field approximation]\label{thm:spherecoat}
Fix $L\geq1$. Let the incident Herglotz density belong to the span of
$\Psi_{lm},X_{lm}$ with $1\leq l\leq L$. There are $h_0,\eta_0>0$
such that, for $0<h\leq h_0$ and $0<\eta<\eta_0$, the
anisotropic vacuum-side coating has a unique scattering solution in
this invariant angular subspace, and
\begin{equation}\label{eq:spherecoatfar}
 \sup_{z\in I}
 \|\FF_{\eta,h}^{c}-\FF_\eta^{\gnor}\|_{L^2(\Sph)}
 \leq C_L(\eta^2+h\eta)\|\mathbf a\|_{L^2(\Sph)}.
\end{equation}
For the isotropic metal-side coating, if $0<h\leq c_*\eta$,
\begin{equation}\label{eq:spheremetalfar}
 \sup_{z\in I}
 \|\FF_{\eta,h}^{m}-\FF_\eta^{\gnor}\|_{L^2(\Sph)}
 \leq C_L(\eta^2+h)\|\mathbf a\|_{L^2(\Sph)}.
\end{equation}
The estimates also hold with $\FF_\eta^{\eff}$ in place of
$\FF_\eta^{\gnor}$. Thus, $h=O(\eta)$ and $h=O(\eta^2)$, respectively,
give first-order far-field equivalence. 
\end{theorem}
\begin{proof}
We give the TM argument first. Let $j_l$ and $h_l^{(1)}$ denote the spherical Bessel and outgoing
spherical Hankel functions, respectively. A regular core state at $r=1$ is
\[
 w_m=\left.
 \binom{rj_l(\kappa\sqrt\eps\,r)}
 {\eps^{-1}\partial_r[rj_l(\kappa\sqrt\eps\,r)]}
 \right|_{r=1},
\]
with a fixed continuous branch of $\sqrt\eps$. Let $w_{\mathrm{reg}}$ and
$w_{\mathrm{out}}$ be the analogous vacuum states with $j_l(\kappa r)$
and $h_l^{(1)}(\kappa r)$, respectively. Let $a_l^{\scat}$ and $a_l^m$ be the scattered and core coefficients.
A unit regular incident mode gives the matching system
\[
 \bigl(w_{\mathrm{out}},-Tw_m\bigr)\binom{a_l^{\scat}}{a_l^m}
 =-w_{\mathrm{reg}}.
\]
At $T=\Id$ its determinant does not vanish. Otherwise a nonzero
regular-core, outgoing-vacuum solution would contradict the energy
identity for $\operatorname{Im}\eps>0$ and Rellich's lemma.
The states are continuous in $z$, and the finite set $1\leq l\leq L$
therefore has uniformly bounded inverse matrices on $I$.
Equations \eqref{eq:spherecoatTerror}--\eqref{eq:coatspherereg} show
that all the matching systems are small perturbations of these local
systems. The inverse identity gives the stated errors in each TM
scattering coefficient.

For TE modes, interchange electric and magnetic roles. With
$\E=uX_{lm}/r$, use state $(u,u')$; its generator is
\[
 A_{l,\mathrm{TE}}(r)=
 \begin{pmatrix}0&1\\\lambda_l/r^2-\kappa^2\eps_\tau(r)&0\end{pmatrix}.
\]
The anisotropic vacuum-side coating is exactly vacuum in this channel,
so its extrapolated transfer is $\Id$. GNOR longitudinal
screening also leaves the TE channel unchanged.
For the isotropic metal-side layer, $\eps_c$ is bounded under
$h\leq c_*\eta$; comparison with the metal generator gives an
$O_L(h)$ transfer error. The same stable matching argument applies.

Finally, the Herglotz map, the mode amplitudes, and the electric far-field
map are bounded between their finite-dimensional angular spaces,
uniformly for $z\in I$. Orthogonality and summation over $(l,m)$ prove
the two estimates and uniqueness within the stated invariant subspace.
\end{proof}

\subsection{A retarded convergence check}

\begin{figure}[htbp]
\centering
\begin{tikzpicture}
\begin{loglogaxis}[width=.47\textwidth,height=5.8cm,
 xlabel={$\eta$},ylabel={coefficient error},
 tick label style={font=\scriptsize},label style={font=\small},
 legend style={font=\scriptsize,at={(0.5,-0.27)},anchor=north,
 draw=none,fill=white},grid=major,
 xmajorgrids=true,ymajorgrids=true,log basis x=10,log basis y=10,
name=leftplot,title={(a) Surface and observation expansions},
 title style={font=\small}]
\addplot+[metalblue,mark=o] coordinates {(1.000000000000e-04,6.509845508364e-08) (1.778279410039e-04,2.059151834329e-07) (3.162277660168e-04,6.514748190713e-07) (5.623413251903e-04,2.061909970768e-06) (1.000000000000e-03,6.530269900044e-06) (1.778279410039e-03,2.070649984071e-05) (3.162277660168e-03,6.579532881764e-05) (5.623413251903e-03,2.098465238068e-04) (1.000000000000e-02,6.737044705599e-04)};
\addlegendentry{first-order surface}
\addplot+[fieldteal,mark=triangle*] coordinates {(1.000000000000e-04,8.756592616041e-08) (1.778279410039e-04,2.769788878991e-07) (3.162277660168e-04,8.762842497754e-07) (5.623413251903e-04,2.773305908021e-06) (1.000000000000e-03,8.782644829202e-06) (1.778279410039e-03,2.784466191384e-05) (3.162277660168e-03,8.845649699029e-05) (5.623413251903e-03,2.820142052219e-04) (1.000000000000e-02,9.048716908121e-04)};
\addlegendentry{first-order observation}
\addplot+[screenorange,mark=square*] coordinates {(1.000000000000e-04,3.683111154299e-11) (1.778279410039e-04,2.071694281922e-10) (3.162277660168e-04,1.165530193152e-09) (5.623413251903e-04,6.559569234885e-09) (1.000000000000e-03,3.694035505499e-08) (1.778279410039e-03,2.082641483832e-07) (3.162277660168e-03,1.176512854494e-06) (5.623413251903e-03,6.670033437388e-06) (1.000000000000e-02,3.805631152999e-05)};
\addlegendentry{second-order observation}
\addplot+[black!65,mark=diamond*,dashed] coordinates {(1.000000000000e-04,1.788761420307e-11) (1.778279410039e-04,1.005962519828e-10) (3.162277660168e-04,5.657643620607e-10) (5.623413251903e-04,3.182229690112e-09) (1.000000000000e-03,1.790204840593e-08) (1.778279410039e-03,1.007418666735e-07) (3.162277660168e-03,5.672376864484e-07) (5.623413251903e-03,3.197264662216e-06) (1.000000000000e-02,1.805754543758e-05)};
\addlegendentry{second-order surface}
\end{loglogaxis}
\begin{loglogaxis}[width=.47\textwidth,height=5.8cm,
 xlabel={$\eta$},ylabel={coefficient error},
 tick label style={font=\scriptsize},label style={font=\small},
 legend style={font=\scriptsize,at={(0.5,-0.27)},anchor=north,
 draw=none,fill=white},grid=major,
 xmajorgrids=true,ymajorgrids=true,log basis x=10,log basis y=10,
at={(leftplot.east)},anchor=west,xshift=1.0cm,
title={(b) Local coating representations},
 title style={font=\small}]
\addplot+[metalblue,mark=o] coordinates {(1.000000000000e-04,6.689749559890e-08) (1.778279410039e-04,2.115994100321e-07) (3.162277660168e-04,6.694226993655e-07) (5.623413251903e-04,2.118512888600e-06) (1.000000000000e-03,6.708400463957e-06) (1.778279410039e-03,2.126492469466e-05) (3.162277660168e-03,6.753364554015e-05) (5.623413251903e-03,2.151867523941e-04) (1.000000000000e-02,6.896924967363e-04)};
\addlegendentry{anisotropic, $h=\eta$}
\addplot+[fieldteal,mark=triangle*] coordinates {(1.000000000000e-04,2.764164678280e-08) (1.778279410039e-04,8.744181014691e-08) (3.162277660168e-04,2.766910863598e-07) (5.623413251903e-04,8.759634990102e-07) (1.000000000000e-03,2.775612327198e-06) (1.778279410039e-03,8.808677259686e-06) (3.162277660168e-03,2.803300858301e-05) (5.623413251903e-03,8.965475646230e-05) (1.000000000000e-02,2.892557070517e-04)};
\addlegendentry{isotropic, $h=0.2\eta^2$}
\addplot+[screenorange,mark=square*] coordinates {(1.000000000000e-04,3.449788809374e-05) (1.778279410039e-04,6.139556236335e-05) (3.162277660168e-04,1.093325213024e-04) (5.623413251903e-04,1.949116324150e-04) (1.000000000000e-03,3.481539815279e-04) (1.778279410039e-03,6.240278695874e-04) (3.162277660168e-03,1.125355327020e-03) (5.623413251903e-03,2.051416571021e-03) (1.000000000000e-02,3.810793682812e-03)};
\addlegendentry{isotropic, $h=0.2\eta$}
\end{loglogaxis}
\end{tikzpicture}
\caption{Retarded TM dipole verification for a unit sphere, with
$z=0.55$, $\gamma=0.08$, and $\sigma=2$, so $\kappa=1.1$.
The vertical axis is the absolute error in the scattered coefficient
for unit regular incident coefficient, relative to the exact
coefficient at $\diff=0$. In (a), the first-order surface model and
first-order observation expansion have slope $2$; the second-order
surface model and second-order observation expansion have slope $3$.
In (b), the anisotropic layer with $h=\eta$ and isotropic layer with
$h=0.2\eta^2$ have slope $2$; the isotropic layer with $h=0.2\eta$ has
slope $1$. These are fixed-mode checks, not tests of uniformity over
all multipoles.}
\label{fig:convergence}
\end{figure}
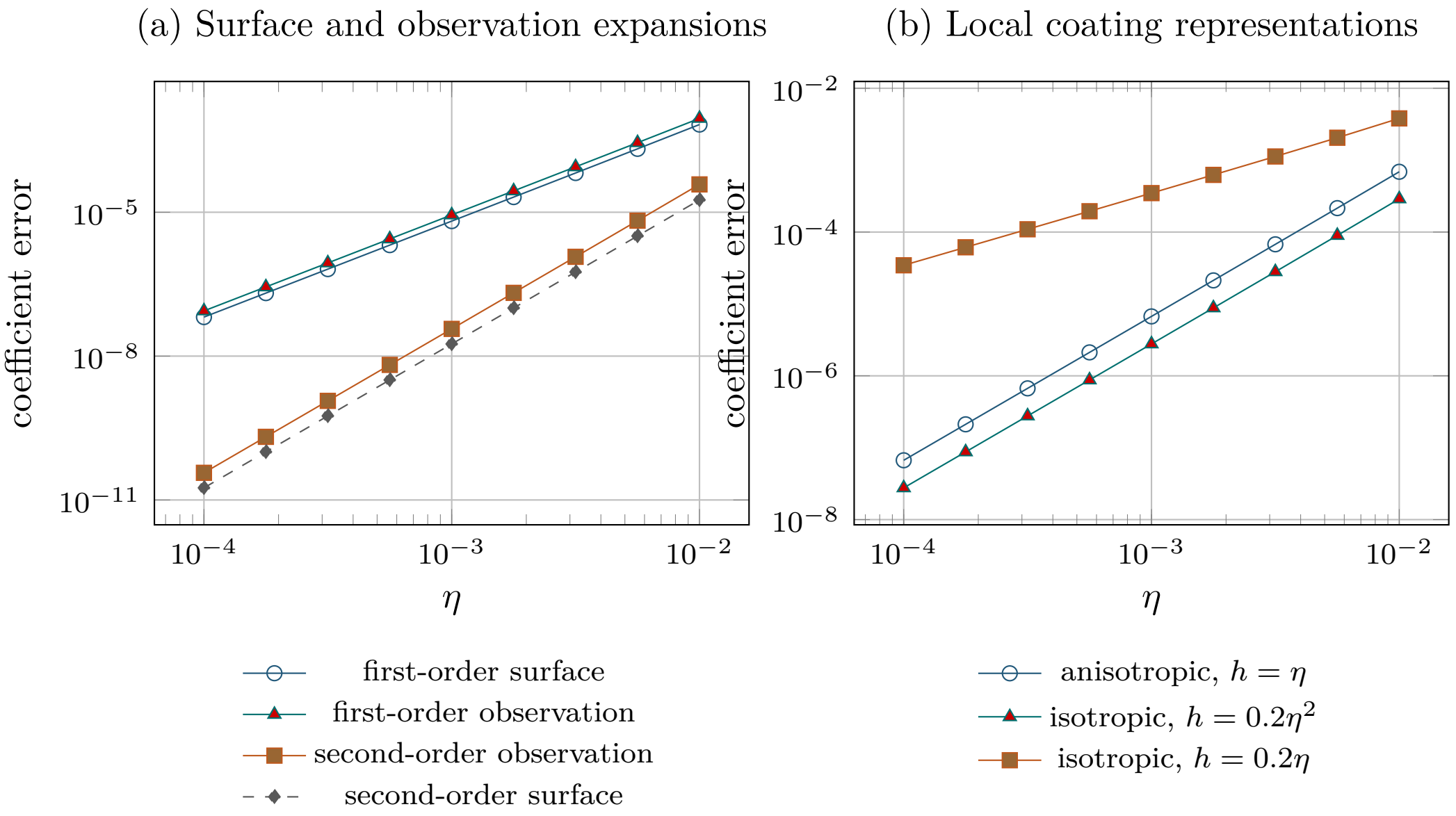

Figure~\ref{fig:convergence} tests the TM dipole coefficient for a unit
sphere at $\kappa=1.1$ and $\diff=0$. The exact and surface coefficients are
obtained from the two-by-two matching system above. The coated
coefficients use the radial system \eqref{eq:coatradialode}, integrated
after rescaling the shell to a unit interval, followed by the stated
reference-surface extrapolation. Thus, the calculation does not resolve
a shrinking radial interval with a fixed physical mesh.

The convergence checks in Figure~\ref{fig:convergence} and Table~\ref{tab:gnorconvergence} use the nine logarithmically spaced
values
\[
 \eta_j=10^{-4+j/4},\qquad j=0,\ldots,8.
\]
The reported convergence rates are least-squares slopes of the logarithm
of the coefficient error against $\log\eta_j$ for $j=0,\ldots,4$.

For clarity, let $a_l(b)$ denote $a_l^{\scat}$, the scattered coefficient for the transfer
$\Id+bB_l(1)$ and write
$a_l(b)=a_l(0)+b a_l^{[1]}+b^2a_l^{[2]}+O(b^3)$.
At $\diff=0$ on the unit sphere, Lemma~\ref{lem:ntdexp} gives
\[
 b_l=\frac\eta{s^3}+\frac{\eta^2}{s^4}+O_l(\eta^3).
\]
The two observation approximations in the figure are therefore
\[
 a_l(0)+\frac\eta{s^3}a_l^{[1]},\qquad
 a_l(0)+\frac\eta{s^3}a_l^{[1]}
       +\eta^2\left(\frac{a_l^{[1]}}{s^4}
                     +\frac{a_l^{[2]}}{s^6}\right).
\]
For the second-order surface model, put $r_l=(1+\eta\lambda_l)^{-1}$.
Since $\mathcal H=1$ on the unit sphere, its exact transfer is
$\Id+\widehat b_l^{(2)}B_l(1)$, where, by \eqref{eq:gnorsphereN},
\[
 \widehat b_l^{(2)}
 =(1-\eps^{-1})\widehat\nu_l^{(2)}
 =(1-\eps^{-1})\left[
     \frac\eta{\sg}(2r_l-r_l^2)+\frac{\eta^2}{\sg^2}r_l^2\right].
\]
This coefficient is inserted in the same matching system; no
observation postprocessing is used for this surface model.
The Taylor coefficients are computed by differentiating the matching
matrix, not by fitting the exact data. A DOP853 integration with
relative tolerance $3\times10^{-13}$ and absolute tolerance
$3\times10^{-15}$ is used for the coating transfers. The fitted slopes
on the five smallest values of $\eta$ are between $2.001$ and $2.002$
for the predicted quadratic errors, $3.001$ for the second-order
observation remainder, $3.000$ for the second-order surface model,
and $1.004$ for the thicker isotropic coating.
These checks support the rates for a fixed retarded mode; the proofs
above provide the analytical estimates. The accompanying script
\texttt{verify\_spherical\_convergence.py} reproduces the data,
fitted slopes, and figure.

The same calculation at $\diff=1$ and $\diff=4$ uses $\Qg$ in the
exact eigenvalue, \eqref{eq:gnorsphereN} for the two surface models,
and $b_0=\eta/(s^2\sg)$ from \eqref{eq:coatb} for the coatings. The
transverse matching states do not depend on $\diff$.
Table~\ref{tab:gnorconvergence} gives the measured slopes, with
$\diff=0$ repeating the corresponding entries of
Figure~\ref{fig:convergence}.
The script \texttt{verify\_gnor\_convergence.py} reproduces the table
and also compares quasistatic poles and residues using $65$-digit
arithmetic. For $(l,\diff)=(1,1)$ and $(2,4)$ with $\gamma=0.08$,
the fitted pole and residue slopes agree with the predicted orders
$2$ and $3$ within $0.01$. A separate collisionless calculation
checks the complex shift in \eqref{eq:gnorzeroloss}.
These finite-mode calculations check the formulas and their numerical
implementation; they do not establish uniformity over modes or parameters.

\begin{table}[htbp]
\centering\small
\begin{tabular}{@{}ccccc@{}}
\toprule
$\diff$ & First-order surface & Second-order surface & Anisotropic & Isotropic\\
\midrule
0 & 2.0013 & 3.0003 & 2.0012 & 2.0017\\
1 & 2.0010 & 3.0000 & 2.0008 & 2.0014\\
4 & 2.0006 & 2.9996 & 2.0004 & 2.0010\\
\bottomrule
\end{tabular}
\caption{GNOR convergence slopes for the retarded TM dipole at
$z=0.55$, $\gamma=0.08$, and $\sigma=2$. The coating thicknesses
are $h=\eta$ and $h=0.2\eta^2$, respectively. Slopes use the five
smallest values in the same $\eta$ grid as Figure~\ref{fig:convergence}.
The predicted orders are $2$, $3$, $2$, and $2$.}
\label{tab:gnorconvergence}
\end{table}

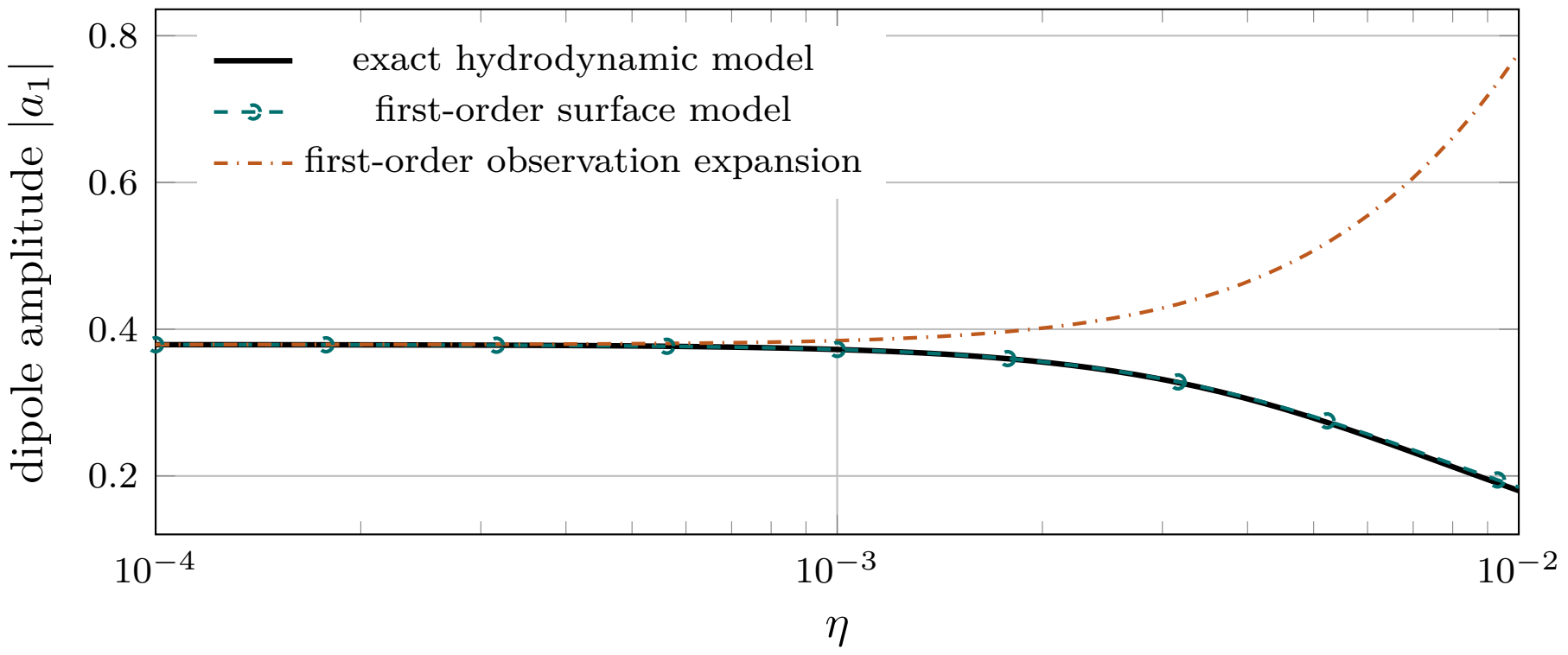
\begin{figure}[htbp]
\centering
\begin{tikzpicture}
\pgfplotstableread[col sep=comma]{
eta,exact,surface,observation
1.000000000000e-04,3.789204607100e-01,3.789205217809e-01,3.790429880763e-01
1.074607828321e-04,3.789055398301e-01,3.789056137379e-01,3.790470292201e-01
1.154781984689e-04,3.788887469037e-01,3.788888364491e-01,3.790521322455e-01
1.240937760752e-04,3.788698250561e-01,3.788699336692e-01,3.790584940348e-01
1.333521432163e-04,3.788484801876e-01,3.788486120714e-01,3.790663443546e-01
1.433012570237e-04,3.788243754270e-01,3.788245357357e-01,3.790759511199e-01
1.539926526059e-04,3.787971247518e-01,3.787973198099e-01,3.790876264829e-01
1.654817099943e-04,3.787662856507e-01,3.787665232229e-01,3.791017338752e-01
1.778279410039e-04,3.787313506853e-01,3.787316403098e-01,3.791186961508e-01
1.910952974970e-04,3.786917377893e-01,3.786920911889e-01,3.791390049984e-01
2.053525026457e-04,3.786467791196e-01,3.786472107082e-01,3.791632318198e-01
2.206734069085e-04,3.785957082475e-01,3.785962357554e-01,3.791920402980e-01
2.371373705662e-04,3.785376454506e-01,3.785382906954e-01,3.792262009143e-01
2.548296747979e-04,3.784715808328e-01,3.784723706704e-01,3.792666077114e-01
2.738419634264e-04,3.783963549655e-01,3.783973224622e-01,3.793142976428e-01
2.942727176209e-04,3.783106367033e-01,3.783118225813e-01,3.793704729016e-01
3.162277660168e-04,3.782128977871e-01,3.782143522053e-01,3.794365266747e-01
3.398208328943e-04,3.781013838022e-01,3.781031685503e-01,3.795140728374e-01
3.651741272548e-04,3.779740810128e-01,3.779762722113e-01,3.796049801730e-01
3.924189758485e-04,3.778286785462e-01,3.778313699699e-01,3.797114117859e-01
4.216965034286e-04,3.776625253575e-01,3.776658325200e-01,3.798358704673e-01
4.531583637601e-04,3.774725813592e-01,3.774766465295e-01,3.799812508767e-01
4.869675251659e-04,3.772553620689e-01,3.772603604257e-01,3.801508995141e-01
5.232991146815e-04,3.770068761069e-01,3.770130232786e-01,3.803486835819e-01
5.623413251903e-04,3.767225548765e-01,3.767301161628e-01,3.805790699707e-01
6.042963902381e-04,3.763971737865e-01,3.764064754188e-01,3.808472157481e-01
6.493816315762e-04,3.760247644564e-01,3.760362073141e-01,3.811590716787e-01
6.978305848599e-04,3.755985174724e-01,3.756125937508e-01,3.815215004610e-01
7.498942093325e-04,3.751106754837e-01,3.751279888872e-01,3.819424115243e-01
8.058421877615e-04,3.745524167520e-01,3.745737068727e-01,3.824309143706e-01
8.659643233601e-04,3.739137297341e-01,3.739399013628e-01,3.829974925833e-01
9.305720409297e-04,3.731832799218e-01,3.732154381208e-01,3.836542007179e-01
1.000000000000e-03,3.723482710360e-01,3.723877628695e-01,3.844148863444e-01
1.074607828321e-03,3.713943038113e-01,3.714427676672e-01,3.852954394881e-01
1.154781984689e-03,3.703052370686e-01,3.703646604962e-01,3.863140715985e-01
1.240937760752e-03,3.690630575985e-01,3.691358445010e-01,3.874916259227e-01
1.333521432163e-03,3.676477675897e-01,3.677368154142e-01,3.888519207359e-01
1.433012570237e-03,3.660373009333e-01,3.661460881583e-01,3.904221262438e-01
1.539926526059e-03,3.642074826645e-01,3.643401663344e-01,3.922331750702e-01
1.654817099943e-03,3.621320489259e-01,3.622935711882e-01,3.943202050385e-01
1.778279410039e-03,3.597827479246e-01,3.599789494191e-01,3.967230314216e-01
1.910952974970e-03,3.571295450136e-01,3.573672815253e-01,3.994866439428e-01
2.053525026457e-03,3.541409567122e-01,3.544282137153e-01,4.026617216044e-01
2.206734069085e-03,3.507845384596e-01,3.511305360941e-01,4.063051559578e-01
2.371373705662e-03,3.470275482746e-01,3.474428270409e-01,4.104805708569e-01
2.548296747979e-03,3.428378023162e-01,3.433342775626e-01,4.152588242834e-01
2.738419634264e-03,3.381847277076e-01,3.387756991710e-01,4.207184758072e-01
2.942727176209e-03,3.330406023254e-01,3.337407040278e-01,4.269462020635e-01
3.162277660168e-03,3.273819506324e-01,3.282070268734e-01,4.340371427357e-01
3.398208328943e-03,3.211910400763e-01,3.221579356361e-01,4.420951614034e-01
3.651741272548e-03,3.144573963966e-01,3.155836537798e-01,4.512330096076e-01
3.924189758485e-03,3.071792320032e-01,3.084826957214e-01,4.615723887869e-01
4.216965034286e-03,2.993646641085e-01,3.008630012667e-01,4.732439131945e-01
4.531583637601e-03,2.910325935462e-01,2.927427504596e-01,4.863869869927e-01
4.869675251659e-03,2.822131254097e-01,2.841507501770e-01,5.011496194792e-01
5.000000000000e-03,2.788615651926e-01,2.808862523261e-01,5.070051127553e-01
5.232991146815e-03,2.729474408641e-01,2.751263099846e-01,5.176882125339e-01
5.623413251903e-03,2.632870745090e-01,2.657185659798e-01,5.361673624605e-01
6.042963902381e-03,2.532926084898e-01,2.559852630615e-01,5.567597231079e-01
6.493816315762e-03,2.430318549608e-01,2.459910609462e-01,5.796459776079e-01
6.978305848599e-03,2.325776526568e-01,2.358054785438e-01,6.050149620143e-01
7.498942093325e-03,2.220054420687e-01,2.255006268558e-01,6.330639760762e-01
8.058421877615e-03,2.113908007777e-01,2.151488967322e-01,6.639993054972e-01
8.659643233601e-03,2.008071140078e-01,2.048207628419e-01,6.980369678985e-01
9.305720409297e-03,1.903235281503e-01,1.945828413447e-01,7.354036829967e-01
1.000000000000e-02,1.800032932113e-01,1.844963014950e-01,7.763380576497e-01
}\resonantdata
\begin{semilogxaxis}[width=.78\textwidth,height=5.6cm,
 xlabel={$\eta$},ylabel={dipole amplitude $|a_1|$},
 xmin=.0001,xmax=.01,grid=major,
 tick label style={font=\scriptsize},label style={font=\small},
 legend style={font=\scriptsize,at={(0.03,0.97)},anchor=north west,
 draw=none,fill=white}]
\addplot[black,very thick] table[x=eta,y=exact]{\resonantdata};
\addlegendentry{exact hydrodynamic model}
\addplot[fieldteal,dashed,thick,mark=o,mark repeat=8,mark size=1.5pt]
 table[x=eta,y=surface]{\resonantdata};
\addlegendentry{first-order surface model}
\addplot[screenorange,dashdotted,thick]
 table[x=eta,y=observation]{\resonantdata};
\addlegendentry{first-order observation expansion}
\end{semilogxaxis}
\end{tikzpicture}
\caption{Retarded TM dipole on the unit sphere at $z=\operatorname{Re}z_0$, with
$\gamma=0.005$, $\sigma=0.35$, $\diff=0$, and unit regular incident coefficient. The surface model
evaluates $a_1(b_0r_1)$, and the observation expansion evaluates
$a_1(0)+b_0a_1^{[1]}$. The exact coefficient is $a_1(b_1)$.}
\label{fig:resonantcmp}
\end{figure}

We also compare the two first-order approximations close to a local
dipole pole. We take $\diff=0$, $\gamma=0.005$, and $\sigma=0.35$.
The local outgoing matching determinant has a zero at
$z_0\simeq0.568269-0.003903\ii$. Figure~\ref{fig:resonantcmp} varies
$\eta$ at the real frequency $z=\operatorname{Re}z_0$. At
$\eta=0.005$, the absolute errors in the complex coefficient are
$0.00300$ for the surface model and $0.228$ for the observation
expansion. This calculation concerns one retarded mode at fixed loss.

\subsection{Scope and passivity}

The planar estimates and Theorem~\ref{thm:spherecoat} are quantitative
coating realizations of the first-order surface law. The sphere theorem
includes curvature and the self-consistent retarded scattering problem,
but its constants depend on the multipole cutoff. It does not establish
uniqueness of the coated problem on the unrestricted angular space.
A general smooth-particle theorem would require uniform control of the
high-tangential-frequency response and of the coated Maxwell resolvent.
Those estimates are not consequences of moment matching.

Passivity is a separate issue even for the hard-wall coefficients.
For every $\diff$,
\begin{equation}\label{eq:coatnonpassive}
 \operatorname{Im}\eps_n^c
 =-\frac{h\operatorname{Im}b_0}{|h+b_0|^2}.
\end{equation}
At $\diff=0$, $z\in(0,1)$, and $\gamma>0$, $\arg s\in(-\pi/4,0)$, so
$\operatorname{Im}b_0>0$ and $\operatorname{Im}\eps_n^c<0$.
Thus the vacuum-side coating in \eqref{eq:anisocoat}, considered as a
standalone local material with the convention $e^{-\ii\omega t}$, is
not passive at these frequencies. For $\diff>0$ the sign of
$\operatorname{Im}b_0$ is not fixed, and no passivity assertion is made.
The fixed-band and fixed-multipole
results use perturbative invertibility and do not infer stability from
coating passivity. Neither a passive realization at all frequencies nor
equality of all scattering poles is proved by this construction.

\section{Classical transmission and observation estimates}\label{app:transmission}

We give the functional setting and the ellipticity argument used in
Lemma~\ref{lem:classical}. All boundaries in this appendix are fixed
and $C^\infty$. We use the smooth Maxwell layer-potential calculus
and its Hodge decomposition; see \cite[Chapter~5]{Nedelec} and
\cite{Buffa,MaxwellTraces,Vodev}. The role of the complex Drude coefficient is
specified below rather than inferred from a positive-coefficient result.

\subsection{Admittance maps and the transmission equation}

Let $e$ be the tangential electric trace. Define
$\mathcal A_z^-e=\nn\times\HH^-$ by the interior local Drude
Maxwell problem with trace $e$, and define
$\mathcal A_z^+e=\nn\times\HH^+$ by the outgoing vacuum problem
with trace $e$. The same normal $\nn$ is used on both sides.
The interior electric Dirichlet problem is uniquely solvable because
$\operatorname{Im}\eps>0$: the imaginary part of its homogeneous
energy identity gives $\E=0$, and the rotated variational form is
coercive after lifting the boundary data. On the unbounded vacuum component, the outgoing electric Dirichlet
problem is uniquely solvable for real $\kappa>0$. It can be constructed
with a combined-field layer representation; interior eigenvalues of an
auxiliary integral representation do not restrict that solution operator.
Hollow particles may enclose bounded vacuum components. Their electric Dirichlet
maps are meromorphic and can have finite-rank poles at cavity
eigenfrequencies. We first work away from those frequencies and then
remove this auxiliary restriction by the finite-dimensional augmentation
described below. On such a regular frequency interval the admittances
exist on the natural Maxwell trace spaces and depend continuously on $z$.

For precision, let $P_\nabla$, $P_{\mathrm{sol}}$, and $P_{\mathrm{h}}$ be the
$L^2$-orthogonal Hodge projections, and let $H^r_\tau(\Gamma)$ be the
Sobolev space of tangential vector fields. The projections act onto
surface gradients, rotated surface
gradients, and harmonic tangential fields, respectively. The last space
is finite dimensional and includes the topological contribution on
each boundary component. Define the scales
\begin{align*}
 \mathcal E^r&=P_\nabla H^r_\tau(\Gamma)
        \oplus P_{\mathrm{sol}}H^{r+1}_\tau(\Gamma)
        \oplus\ran P_{\mathrm{h}},\\
 \mathcal M^r&=P_\nabla H^{r+1}_\tau(\Gamma)
        \oplus P_{\mathrm{sol}}H^r_\tau(\Gamma)
        \oplus\ran P_{\mathrm{h}}.
\end{align*}
At $r=-1/2$ these are, with equivalent norms,
$H^{-1/2}_{\mathrm{curl}}(\Gamma)$ and
$H^{-1/2}_{\mathrm{div}}(\Gamma)$, respectively. The first notation
means that the tangential trace has surface curl in $H^{-1/2}$;
the second requires surface divergence in $H^{-1/2}$.
We use $\curl_\Gamma e=-\diver_\Gamma(\nn\times e)$.
The natural electric trace is $\E_\tau\in\mathcal E^{-1/2}$,
and the rotated magnetic trace is $\nn\times\HH\in\mathcal M^{-1/2}$;
see \cite{MaxwellTraces}. The admittances extend
as continuous maps $\mathcal A_z^\pm:\mathcal E^r\to\mathcal M^r$
for all real $r$, by the smooth-boundary calculus and duality.

Set $\mathcal B_z=\mathcal A_z^+-\mathcal A_z^-$. Prescribing
$e^+-e^-=\nabla_\Gamma v$ and continuity of the tangential magnetic
field is equivalent to
\begin{equation}\label{eq:admittancetransmission}
 \mathcal B_z e^-=-\mathcal A_z^+\nabla_\Gamma v,
 \qquad e^+=e^-+\nabla_\Gamma v.
\end{equation}
We next check the ellipticity needed to invert this equation.

\subsection{Principal symbols, index, and normal response}

Flatten the interface with the metal in $t<0$, and take a nonzero
tangential covector $\xi$. Put $k=|\xi|$ and, for sufficiently large
$k$, choose the roots
\[
 q_0=(k^2-\kappa^2)^{1/2},\qquad
 q_m=(k^2-\eps\kappa^2)^{1/2},\qquad \operatorname{Re}q_j>0.
\]
The decaying factors are $e^{-q_0t}$ outside and $e^{q_mt}$ inside.
In the basis parallel and perpendicular to $\xi$, the frozen planar
admittances are
\begin{equation}\label{eq:admittancesymbols}
 \mathcal A_z^+\sim
 \begin{pmatrix}\ii\kappa/q_0&0\\0&-\ii q_0/\kappa\end{pmatrix},
 \qquad
 \mathcal A_z^-\sim
 \begin{pmatrix}-\ii\kappa\eps/q_m&0\\0&\ii q_m/\kappa\end{pmatrix}.
\end{equation}
Thus the gradient channel has order $-1$ and the rotated-gradient
channel has order $1$. These are mixed-order symbols: one must use
the Hodge scales, not regard both channels as scalar order-one maps.

Let
$\Lambda=(\Id+\Delta_{\mathrm{H}})^{1/2}$, where $\Delta_{\mathrm{H}}$ is the
nonnegative Hodge Laplacian on tangential one-forms. It commutes with
the Hodge projections. We write $\Psi^m(\Gamma)$ for the classical
pseudodifferential operators of order at most $m$ on tangential fields;
operators with smooth kernels belong to every such class.

\begin{lemma}[Hodge-block estimates]\label{lem:hodgeblocks}
Away from the auxiliary cavity eigenfrequencies, the admittances have
the following blocks, modulo operators with smooth kernels:
\begin{align}
 P_\nabla\mathcal A_z^\pm P_\nabla
   &=a_\pm P_\nabla\Lambda^{-1}P_\nabla
       +\Psi^{-2}(\Gamma),\label{eq:gradientblock}\\
 P_{\mathrm{sol}}\mathcal A_z^\pm P_{\mathrm{sol}}
   &=b_\pm P_{\mathrm{sol}}\Lambda P_{\mathrm{sol}}
       +\Psi^0(\Gamma),\label{eq:solenoidalblock}\\
 P_{\mathrm{sol}}\mathcal A_z^\pm P_\nabla,\quad
 P_\nabla\mathcal A_z^\pm P_{\mathrm{sol}}
   &\in\Psi^{-1}(\Gamma),\label{eq:offdiagonalblocks}
\end{align}
where $a_+=\ii\kappa$, $a_-=-\ii\kappa\eps$,
$b_+=-\ii/\kappa$, and $b_-=\ii/\kappa$.
Here an additive operator class denotes a remainder in that class.
The estimates are uniform on compact frequency sets where the
admittances are regular. Blocks involving $P_{\mathrm h}$ have smooth
kernels.
\end{lemma}
\begin{proof}
We first justify the improved order for gradient electric data on a
curved boundary. In either homogeneous medium, write $\eps_j=\eps$
or $1$ for its permittivity. For $e=P_\nabla e$, Maxwell's equations
give
\[
 (-\Delta-\kappa^2\eps_j)\HH=0,\qquad
 H_n=(\ii\kappa)^{-1}\curl_\Gamma e=0,\qquad
 (\curl\HH)_\tau=-\ii\kappa\eps_j e.
\]
Extend $\nn$ along normal lines and put
$\mathcal S=\nabla_\Gamma\nn$, restricted to tangent vectors.
The boundary identity
\[
 \nn\times\curl\HH
 =\nabla_\Gamma H_n-\partial_n\HH_\tau-\mathcal S\HH_\tau
\]
then gives
\[
 H_n=0,\qquad
 \partial_n\HH_\tau+\mathcal S\HH_\tau
       =\ii\kappa\eps_j\nn\times e.
\]
This is an elliptic mixed boundary problem for the vector Helmholtz
operator: its principal boundary conditions are Dirichlet in the normal
component and Neumann in the tangential components. In boundary normal
coordinates, its decaying principal solutions have normal exponent
$|\xi|$. Thus the tangential Neumann-to-Dirichlet symbol has order
$-1$. The boundary parametrix, obtained by solving successively for
the lower-order symbols, has a smooth remainder; curvature enters the
Robin coefficient $\mathcal S$ and the lower-order coordinate terms.
Consequently $\mathcal A_z^\pm P_\nabla\in\Psi^{-1}(\Gamma)$.
Its leading symbol is obtained by freezing these coefficients and is
$a_\pm|\xi|^{-1}$ on the gradient channel, as in
\eqref{eq:admittancesymbols}. The next symbol has order at most $-2$,
which proves \eqref{eq:gradientblock} and the first off-diagonal bound.
This argument uses the local elliptic parametrix, not invertibility of
an auxiliary vector Helmholtz problem; any global finite-dimensional
corrections have smooth kernels.

For arbitrary electric data, $H_n=(\ii\kappa)^{-1}\curl_\Gamma e$
has order one in $e$. The same mixed problem therefore gives an
admittance of order at most one. Its leading symbol on the
rotated-gradient channel is $b_\pm|\xi|$, and subtraction leaves
order at most zero, proving \eqref{eq:solenoidalblock}.
Lorentz reciprocity gives
\[
 \int_\Gamma e_1\cdot\mathcal A_z^\pm e_2\dd S
 =\int_\Gamma e_2\cdot\mathcal A_z^\pm e_1\dd S
\]
for smooth data. This is the bilinear transpose identity, without
complex conjugation; for the exterior problem the boundary term at
infinity vanishes for two outgoing fields. The Hodge projections are
real and symmetric under the same pairing. Transposing the first
off-diagonal block thus proves the second bound in
\eqref{eq:offdiagonalblocks}. Finally, harmonic fields are smooth and
finite dimensional, and interactions between disjoint boundary
components have smooth kernels. Continuity of the local symbols and
of the global solution operators gives the asserted uniformity.
\end{proof}

We now perform the order reduction. The isomorphisms
\[
 J_E=P_\nabla\Lambda^{-r}
       +P_{\mathrm{sol}}\Lambda^{-r-1}+P_{\mathrm{h}},\qquad
 J_M=P_\nabla\Lambda^{-r-1}
       +P_{\mathrm{sol}}\Lambda^{-r}+P_{\mathrm{h}}
\]
map $L^2_\tau(\Gamma)$ onto $\mathcal E^r$ and $\mathcal M^r$.
Lemma~\ref{lem:hodgeblocks} gives
\begin{equation}\label{eq:reducedsymbol}
 J_M^{-1}\mathcal B_zJ_E
 =\ii\kappa(1+\eps)P_\nabla
       -\frac{2\ii}{\kappa}P_{\mathrm{sol}}+P_{\mathrm{h}}+\mathcal K_z,
\end{equation}
where $\mathcal K_z$ is compact on $L^2_\tau$ and continuous in $z$.
Indeed, after multiplication by $J_M^{-1}$ and $J_E$, the diagonal
remainders in \eqref{eq:gradientblock}--\eqref{eq:solenoidalblock}
have order at most $-1$. The two off-diagonal blocks retain order at
most $-1$. They are therefore compact on $L^2_\tau$.
Replacing the harmonic block by $P_{\mathrm h}$ changes only a
finite-rank operator. This proves \eqref{eq:reducedsymbol} on the
curved boundary, with the different channel orders retained.

The explicit operator preceding $\mathcal K_z$ is invertible because
$\kappa\ne0$ and $1+\eps\ne0$. Hence $\mathcal B_z$ is Fredholm
of index zero. If $\mathcal B_ze=0$, the two fields with common trace
$e$ form an outgoing homogeneous transmission solution. The outgoing
energy identity and $\operatorname{Im}\eps>0$ force the interior
electric field and the exterior far field to vanish. Rellich's lemma gives the zero field on the unbounded vacuum
component. On a bounded vacuum component, both tangential Cauchy traces
vanish, so boundary uniqueness gives the zero field there as well.
Hence $e=0$. Elliptic regularity
makes any distributional null vector smooth, so this argument applies
on every Hodge scale. Thus,
\[
 \mathcal B_z^{-1}:\mathcal M^r\longrightarrow\mathcal E^r
\]
is bounded for every real $r$. This also proves existence for \eqref{eq:admittancetransmission}
at the regular frequencies.

To include a cavity eigenfrequency $z_c\in I$, let $G$ be a bounded
vacuum component, with outward normal $\boldsymbol\nu_G=-\nn$. Let
$\mathscr M_G$ be the nonnegative self-adjoint realization of
$\curl\curl$ in $L^2(G)^3$ with zero tangential electric trace.
Its positive eigenvalues are discrete with finite multiplicity; its
zero eigenspace causes no difficulty because $\kappa$ stays away from
zero. Suppose $\lambda_c=\sigma^2z_c^2$ is an eigenvalue, and let
$\boldsymbol\phi_1,\ldots,\boldsymbol\phi_q$ be an orthonormal basis of its eigenspace.
These eigenfields are smooth up to the boundary. Write $P_c$ for the
orthogonal projection onto their span and
$d_c(z)=\lambda_c-\kappa(z)^2$.

For smooth tangential data $e$ on $\partial G$, choose a smooth lifting
$\mathcal R_Ge$ with tangential trace $e$ and $P_c\mathcal R_Ge=0$;
subtracting $P_c$ from any lifting preserves its boundary trace.
Define
\begin{align*}
 F_{G,z}e&=-(\curl\curl-\kappa^2)\mathcal R_Ge,\\
 \E_{G,z}^{\mathrm{reg}}e
 &=\mathcal R_Ge+
   (\mathscr M_G-\kappa^2)^{-1}_{P_c^\perp}
       (\Id-P_c)F_{G,z}e.
\end{align*}
The inverse on $P_c^\perp$ is regular in a neighborhood of $z_c$
containing no other positive cavity eigenfrequency. The field
$\E_{G,z}^{\mathrm{reg}}e$ is orthogonal to $\ran P_c$ and satisfies
\[
 (\curl\curl-\kappa^2)\E_{G,z}^{\mathrm{reg}}e
 =-\sum_{j=1}^q(\mathcal D_Ge)_j\boldsymbol\phi_j,
 \qquad
 (\mathcal D_Ge)_j=(F_{G,z}e,\boldsymbol\phi_j)_{L^2(G)}.
\]
Green's identity and the orthogonality of the lifting give the
explicit boundary functional
\begin{equation}\label{eq:cavitycompatibility}
 (\mathcal D_Ge)_j
 =-\int_{\partial G}(\boldsymbol\nu_G\times e)
            \cdot\overline{\curl\boldsymbol\phi_j}\dd S.
\end{equation}
In particular, $\mathcal D_G$ is independent of $z$ and is continuous
on every boundary Sobolev scale. With $c=(c_1,\ldots,c_q)\in\C^q$,
the complete cavity electric field is
\[
 \E_G=\E_{G,z}^{\mathrm{reg}}e+\sum_{j=1}^q c_j\boldsymbol\phi_j,
 \qquad d_c(z)c=\mathcal D_Ge.
\]
This equation both reconstructs the field away from $z_c$ and imposes
the compatibility condition at $z_c$.

Replace the resonant cavity contribution to $\mathcal A_z^+$ by
$\nn\times\curl\E_{G,z}^{\mathrm{reg}}/(\ii\kappa)$, keeping
the usual admittances on the other vacuum components, and call the
result $\mathcal A_z^{+,\mathrm{reg}}$. Set
\[
 \mathcal B_z^{\mathrm{reg}}
    =\mathcal A_z^{+,\mathrm{reg}}-\mathcal A_z^-,
 \qquad
 \mathcal U_zc
    =\sum_{j=1}^q c_j
       \frac{\nn\times\curl\boldsymbol\phi_j}{\ii\kappa}
       \quad\hbox{on }\partial G,
\]
and extend $\mathcal U_zc$ by zero to the remaining components of
$\Gamma$. In the following system, $\mathcal D_G$ acts on the
restriction to $\partial G$. For a prescribed scalar jump $v$, the
augmented transmission equation is
\begin{equation}\label{eq:cavityaugmentation}
 \begin{pmatrix}
   \mathcal B_z^{\mathrm{reg}}&\mathcal U_z\\
   -\mathcal D_G&d_c(z)\Id_q
 \end{pmatrix}
 \begin{pmatrix}e^-\\c\end{pmatrix}
 =\begin{pmatrix}
   -\mathcal A_z^{+,\mathrm{reg}}\nabla_\Gamma v\\
   \mathcal D_G\nabla_\Gamma v
 \end{pmatrix}.
\end{equation}
The first row is continuity of the magnetic trace, with
$e^+=e^-+\nabla_\Gamma v$; the second is the cavity compatibility
equation. For $d_c(z)\ne0$, eliminating $c$ gives precisely
\eqref{eq:admittancetransmission}.

The regular part has the same local principal symbols as the original
admittance: the removed poles and the smooth finite-rank forcing above
do not change them. Thus $\mathcal B_z^{\mathrm{reg}}$ is Fredholm
of index zero from $\mathcal E^r$ to $\mathcal M^r$.
The operator in \eqref{eq:cavityaugmentation} differs by a finite-rank
operator from $\diag(\mathcal B_z^{\mathrm{reg}},\Id_q)$ and hence
has index zero from $\mathcal E^r\oplus\C^q$ to
$\mathcal M^r\oplus\C^q$.
Every vector in its kernel reconstructs a homogeneous physical
transmission field, including at $z_c$, so the energy and
boundary-uniqueness argument above makes all fields zero.
It follows that $e^-=0$ and then $c=0$, since the regular cavity field
is orthogonal to the independent eigenfields. The augmented operator
is therefore invertible. If several cavities resonate at $z_c$, take
the direct sum of their eigenspaces and compatibility equations.
Only finitely many positive cavity eigenvalues meet the bounded
frequency window. Continuity of these augmented inverses and a finite
cover of $I$ give the required uniform estimates without excluding
vacuum cavities.

At a regular frequency the normal component is recovered from the
identity $\diver_\Gamma(\nn\times\HH^-)=\ii\kappa\eps E_n^-$. Hence
\begin{equation}\label{eq:Tfactorization}
 \TT v=-\frac{\eps-1}{\ii\kappa\eps}
 \diver_\Gamma\mathcal A_z^-
 \mathcal B_z^{-1}\mathcal A_z^+\nabla_\Gamma v.
\end{equation}
Since $\nabla_\Gamma v\in\mathcal E^r$ for $v\in H^{r+1}$ and
$\diver_\Gamma:\mathcal M^r\to H^r$, this proves the first mapping
property in \eqref{eq:classicalmaps}. In the flattened gradient channel,
\[
 \widehat{\TT v}(\xi)
 =\frac{(1-\eps)k^2}{q_m+\eps q_0}\widehat v(\xi)
\]
at frozen coefficients. Its principal symbol is
$(1-\eps)k/(1+\eps)$, as stated in Lemma~\ref{lem:classical}.
Consequently $\TT$ and its $L^2$ adjoint have order one on the ordinary
scalar Sobolev scale. Constants in $v$ on each component produce no
electric jump and belong to the kernel of $\TT$ and $\LL$.

\subsection{Smoothing observations and incident fields}

An exterior Maxwell representation expresses the field and its far
field as surface potentials of the Cauchy traces
$(e^+,\mathcal A_z^+e^+)$. Their kernels are smooth on
$\overline K\times\Gamma$ for any fixed observation domain $K$ at
positive distance from $\Gamma$, and the radiation kernels are smooth
on $\Sph\times\Gamma$. The trace reconstruction in
\eqref{eq:admittancetransmission} is continuous on every Sobolev scale.
Let $\mathcal O_{z,K}v$ denote the exterior electric field restricted to
$K$. Composition therefore yields, for all integers $m,N\geq0$,
\[
 \|\LL v\|_{L^2(\Sph)}+
 \|\mathcal O_{z,K}v\|_{H^m(K)}
 \leq C_{m,N,K}\|v\|_{H^{-N}(\Gamma)}.
\]
These maps on distributions are the unique continuous extensions
from smooth $v$. In particular, $\LL^*$ maps $L_t^2(\Sph)$ into
$H^N(\Gamma)$ for every $N$. The same statement holds for the
adjoint of $\mathcal O_{z,K}$ with its appropriate Hilbert-space dual.

For the incident problem, separate the known incident Cauchy data and
apply the same inverse transmission operator to the remaining jumps.
Every fixed-order derivative of a Herglotz wave on a bounded set is
bounded by $C\|\mathbf a\|_{L^2(\Sph)}$, uniformly on $I$.
Smooth elliptic regularity and the normal trace formula then give
$\|f_0\|_{H^r}\leq C_r\|\mathbf a\|$ for every $r$.
At an auxiliary cavity eigenfrequency the same reconstruction uses
the augmented inverse; its finite-dimensional blocks are smoothing.
Thus neither the scalar principal symbol nor the observation estimates
change. This proves all the assertions of Lemma~\ref{lem:classical}.

\end{document}